\pdfoutput=1
\documentclass[12pt]{article}
\usepackage{amsmath,amsfonts}

\usepackage{fullpage,amsmath} % 页边距变小
\usepackage{graphicx,kantlipsum} % 加图
\usepackage{cite} % 可以将[1,2,3] 变为[1-3]
\usepackage{indentfirst} % 解决章节首段首行不缩进的问题
\usepackage{algorithm} 
\usepackage{algorithmic}  
\usepackage[algo2e]{algorithm2e} 

\usepackage{authblk} 
\usepackage{bm} % \bm加粗

\usepackage[outdir=./]{epstopdf}
\usepackage{float}
\usepackage{tikz-cd}

\DeclareGraphicsExtensions{.eps}

\date{}                     %% if you don't need date to appear

\author[1,2,a]{Xiao-Kai An}
\author[1,2]{Lin Du}
\author[1,2,3]{Zi-Chen Deng}
\author[2]{Yu-jia Zhang}
\affil[1]{MIIT Key Laboratory of Dynamics and Control of Complex Systems, Northwestern Polytechnical University, Xi'an 710072, China}
\affil[2]{School of Mathematics and Statistics, Northwestern Polytechnical University, Xi'an 710072, China}
\affil[3]{School of Aeronautics, Northwestern Polytechnical University, Xi'an 710072, China}

\affil[a]{anxiaokai@mail.nwpu.edu.cn}
\title{Residual-Based Multi-peak Sampling Algorithm in Inverse Problems of Dynamical Systems}
\begin{document}

\maketitle

\begin{abstract}
% Stochastic differential equations can describe a wide range of dynamical systems, and obtaining the governing equations of these systems is a prerequisite for studying the nonlinear dynamical behavior of the system. In 2022, Jared O'Leary et al. used a neural network to fit the propagation of system moment information to obtain the equation of the system. This method alleviates the computational difficulty of traditional methods to a certain extent, but still relies on a large amount of trajectory data. Aiming at this shortcoming, we propose a residual-based multi-peak sampling algorithm. Evaluate the training results of each generation of the neural network, design and calculate the residual according to the noise intensity, perform multi-peak sampling where the residual is large, and add the samples to the training set to retrain the neural network. In order to prevent the neural network from falling into the trap of overfitting, we discretize the sampling points. We conduct case studies using first- and second-order nonlinear dynamical systems, and perform bifurcation and first crossing probability analyzes of the fitted equations. The results show that our proposed sampling strategy can reconstruct the stochastic dynamical behavior of the system with only a small number of sample points. Finally, we test the algorithm by adding disturbance noise in the data, and the results show that the sampling strategy has better numerical robustness and stability.

Stochastic differential equations can describe a wide range of dynamical systems, and obtaining the governing equations of these systems is the premise of studying the nonlinear dynamic behavior of the system. Neural networks are currently the most popular approach in the inverse problem of dynamical systems. In order to obtain accurate dynamical equations, neural networks need a large amount of trajectory data as a training set. To address this shortcoming, we propose a residual-based multi-peaks sampling algorithm. Evaluate the training results of each epoch of neural network, calculate the probability density function $P(x)$ of the residual, perform sampling where the $P(x)$ is large, and add samples to the training set to retrain the neural network. In order to prevent the neural network from falling into the trap of overfitting, we discretize the sampling points. We conduct case studies using two classical nonlinear dynamical systems and perform bifurcation and first escape probability analyzes of the fitted equations. Results show that our proposed sampling strategy requires only 20$\sim $30\% of the sample points of the original method to reconstruct the stochastic dynamical behavior of the system. Finally, the algorithm is tested by adding interference noise to the data, and the results show that the sampling strategy has better numerical robustness and stability.

\end{abstract}

% \begin{enabstract}
% English abstract

% %“\par在段首，表示另起一行，“\textbf{}”,花括号内的内容加粗显示
% \end{enabstract}

% A helpful LaTeX table generator: \\
% \url{https://www.tablesgenerator.com/} 

\newpage
\section{Introduction}

Because of the rich dynamic behavior of nonlinear systems, more and more scholars use nonlinear systems instead of linear systems in mathematical modeling to meet more stringent requirements\cite{boechler2011bifurcation,antonio2012frequency,strachan2013subharmonic,kecik2014parametric,vakakis2008nonlinear,green2012benefits,amin2012powering,quinn2011energy}. Usually it is difficult to obtain the nonlinear dynamic equation of a system only by mechanism. In the 1990s, the concept of hybrid modeling was proposed, and scholars tried to integrate data and mechanisms to model unknown systems\cite{psichogios1992hybrid,rico1992discrete,thompson1994modeling}. With the rapid development of data science, data-driven modeling has become Powerful alternative and complementary methods for first-principles modeling\cite{liu2022complex}. Over the past decade, neural network-based recognition has clearly been the most popular black-box modeling technique in the field of dynamics\cite{noel2017nonlinear}. Machine learning methods including neural networks have been gradually applied to the inverse problems of ordinary differential equations, partial differential equations and stochastic differential equations\cite{pei2015demonstration,ayala2016cascaded,tavakolpour2015parametric,paduart2010identification,dreesen2015decoupling}.

Usually, the use of neural networks to fit dynamic equations will face three challenges: (1) the number of parameters of neural networks is large; (2) the training data of neural networks is difficult to obtain; (3) the interpretability and robustness of neural network training results Poor stickiness. Using neural networks embedded with physical information will significantly alleviate the above three problems, such as Hamiltonian neural networks\cite{greydanus2019hamiltonian}, Lagrangian neural networks \cite{cranmer2020lagrangian}, Physics-informed neural networks\cite{raissi2019physics}, Neural ODEs\cite{chen2018neural}, etc. 

Despite the good performance of these neural networks, we note that there is still some room for improvement in designing sampling algorithms. A stochastic physics-informed neural ordinary differential equation framework (SPINODE) was proposed recently, which can effectively fit stochastic physics-informed neural ordinary differential equations with Gaussian noise excitation\cite{o2022stochastic}. In the fitting task using SPINODE, due to the existence of random terms in the equation, we need a large amount of trajectory data from the same initial value to obtain the statistical characteristics of the system. Statistical data from each sample point is expensive, so it is very valuable to use an efficient sampling algorithm that uses fewer samples while maintaining the accuracy of the fit. In practical engineering, this means that fewer sensors can be used to avoid damage to the device structure and save costs. Therefore, the main task of this paper is to design an efficient sampling algorithm suitable for SPINODE to achieve higher fitting accuracy with less sample size.

In the sampling step of the inverse problem of the dynamical system, the Residual-Based Adaptive Refinement (RAR) algorithm is the most simple and efficient method undoubtedly\cite{lu2021deepxde}. The RAR algorithm emphasizes the selection of $k$ points with the largest residual error for sampling, which is considered to be the original residual-based sampling algorithm. Based on this description, many scholars have proposed more sophisticated residual sampling algorithms recently. These new methods can be roughly divided into the following two areas:

\begin{enumerate}
  \item Method for computing the probability density function $P(x)$ of the residual:
  
The definition of $P(x)$, which is proportional to the error, was first proposed, as in Eq.\eqref{p_x:1}\cite{nabian2021efficient}. However, this definition only focuses on the position with a large residual error, and it is easy to ignore the position with a small residual error. In order to extract more information from the residual, Residual-based adaptive distribution (RAD) algorithm was proposed\cite{wu2023comprehensive}. The RAD algorithm uses Eq.\eqref{p_x:2} to calculate the $P(x)$ of the residual ($k$ and $c$ are two hyperparameters), which is a highly generalized formula. The methods mentioned in Ref. \cite{nabian2021efficient,gao2023active,tang2023pinns} are all special cases of the RAD algorithm. In addition, a definition similar to RAD in the form of Eq.\eqref{p_x:3} has been proposed\cite{hanna2022residual}.

% Mohammad Amin Nabian et al proposed to use formula 1 to calculate residual Poor probability distribution function P(x)\cite{nabian2021efficient}. On this basis, Chenxi Wu et al. proposed the Residual-based adaptive distribution (RAD) algorithm\cite{wu2023comprehensive}. The RAD algorithm uses formula 2 to calculate the probability distribution function P(x) of the residual (k and c are two hyperparameters), which is a highly generalized formula. The methods mentioned in literature \cite{nabian2021efficient},\cite{gao2023active} and \cite{tang2023pinns} are all RAD A special case of the algorithm. In addition, John M Hanna et al. proposed to use Equation 3 to calculate the probability distribution function P(x) of the residual\cite{hanna2022residual}.

\begin{subequations} \label{p_x}
\begin{equation} P(x)\propto \varepsilon (x) \label{p_x:1}\end{equation}
\begin{equation} p(\mathbf{x}) \propto \frac{\varepsilon^{k}(\mathbf{x})}{\mathbb{E}\left[\varepsilon^{k}(\mathbf{x})\right]}+c \label{p_x:2} \end{equation}
\begin{equation}
    P(x)\propto max\left \{  log(\varepsilon(x)/\varepsilon_{0}   ),0\right \} \label{p_x:3}
\end{equation}
\end{subequations}

  \item Sampling method based on P(x):

In order to disperse the sampling locations as much as possible, the sampling space can be divided into multiple sub-regions and then sampled on each sub-region\cite{zeng2022adaptive}. Node generation techniques are applied to obtain more samples in regions with larger residuals\cite{peng2022rang}. Removing points with small residuals and adding points with large residuals is a good sampling strategy\cite{zapf2022investigating}.

\end{enumerate}

No matter what formula is used to calculate the probability distribution function $P(x)$ of the residual, it is necessary to face the specific problem of how to sample in the large residual area. Selecting the $k$ points with the largest residual error for sampling is in the RAR algorithm, but this method lacks details and does not consider the impact of multiple maximum values. Segmenting the entire sampling space into several subregions for sampling is an improvement, but it may be necessary to divide regions of different sizes in different tasks, and the hyperparameters lead to low generalization ability of the method. To improve the efficiency and generalization performance of sampling methods, we design two new residual-based multimodal sampling algorithms (RBMS-I and RBMS-II). The RBMS algorithm is inspired by the peak filter. The RBMS-II algorithm is easy to implement and has no hyperparameters, it is a fully adaptive algorithm. We combine the RBMS algorithm and the SPINOD0E algorithm to identify hidden physics in stochastic differential equations, using a grazing system and an RVDP system for case studies. In order to further verify the effectiveness of the RBMS algorithm, we analyzed the stochastic dynamic behavior of the system using neural network parametric equations.

The paper is organized as follows: Section 2 describes the process of combining RBMS algorithm with SPINOD0E algorithm to identify hidden physics in stochastic differential equations. Two cases are studied using the RBMS algorithm in Section 3: a grazing system (first order), and an RVDP system (second order). Section 4 discusses the robustness of the RBMS algorithm in the face of initial disturbances. Concluding comments are made in Section 5.

\section{Method}
A dynamical system excited by Gaussian white noise can be defined by the following equations of motion in state space.
\begin{equation}
    \frac{\mathrm{d}}{\mathrm{d} t} X_j=f_j(\boldsymbol{X})+\sum_{l=1}^m g_{j l}(\boldsymbol{X}) W_l(t), \quad j=1,2, \cdots, n, \label{eq1}
\end{equation}

Where $\boldsymbol{X}=\left[X_1, X_2, \cdots, X_n\right]^{\mathrm{T}}, W_j(t)$ is Gaussian white noise, and its correlation function is
\begin{equation}
    E\left[W_l(t) W_s(t+\tau)\right]=2 \pi K_{l s} \delta(\tau) \quad l, s=1,2, \cdots, m . \label{eq2}
\end{equation}
For ease of analysis we convert the form of the system shown in \eqref{eq1} to the It$ \hat{\rm o}$  form.
\begin{equation}
\mathrm{d} X(t)=f(X)dt+g(X)d W(t) \label{eq3}
\end{equation}

\subsection{Data collection}
% Acquiring large and realistic dynamical system data is usually expensive, so we use stochastic simulation methods to obtain trajectory data of known systems. 
% For stochastic differential equations, commonly used numerical simulation solutions include Euler method, Milstein method, stochastic Longo-Kutta method, etc. Among them, the Euler method is 0.5-order convergent, the Milstein method is 1-order weakly convergent, and the stochastic Runge-Kutta method is 1-order strongly convergent. High-precision numerical methods usually take more time. The problem studied in this paper involves a large amount of data, so we use the Euler-Milstein method shown in formula (1) to perform numerical simulation of stochastic differential equations, taking into account the accuracy requirements and time costs.

The problem studied in this paper requires a large amount of trajectory data, and we use the Euler-Milstein method shown in Eq.\eqref{5} to numerically simulate the stochastic differential equation. This method has advantages in solution accuracy and time cost.
% For each $X$ we can use Equation \eqref{4} to obtain a numerical solution for the number of $L$ steps.
\begin{equation}
    X_j=X_{j-1}+f\left(X_{j-1}\right) \Delta t+g\left(X_{j-1}\right)\left(W\left(\tau_j\right)-W\left(\tau_{j-1}\right)\right), j=1,2, \ldots, L  \label{5}
\end{equation}

% Although starting from the same initial value, the trajectories of numerical solutions of stochastic differential equations are different because of the existence of diffusion terms. It is therefore necessary to look for regularities by means of statistics, and here we calculate the first and second moments of these trajectories. Higher order moments can be used in future studies.The calculation methods of first-order moment and second-order moment are shown in formulas \eqref{eq5} and \eqref{eq6}.
Despite starting from the same initial values, the trajectories of numerical solutions of stochastic differential equations are different due to the presence of the diffusion term. Therefore we use the first and second moments of the trajectory data as statistics. The calculation methods of first-order moment and second-order moment are shown in Eq.\eqref{eq6} and Eq.\eqref{eq7}.

\begin{equation}
\mu_x\left(t_k\right)=\frac{1}{N} \sum_{n=1}^N x_n\left(t_k\right), \label{eq6}
\end{equation}
\begin{equation}
    \Sigma_x\left(t_k\right)=\frac{1}{N} \sum_{n=1}^N\left(x_n\left(t_k\right)-\mu_x\left(t_k\right)\right)\left(x_n\left(t_k\right)-\mu_x\left(t_k\right)\right)^{\top} \label{eq7}
\end{equation}

Where $n \in\{1, \ldots, N\}$ denotes the trajectory index. We divide the initial value space with a grid, and select a set of initial values at each node of the grid to perform $N$ simulations, and the step size of each simulation is $L$.
% Where $n \in\{1, \ldots, N\}$ denotes the trajectory index.In order to provide enough training data for the neural network, we use a grid for sampling in the sampling interval. The initial value for each simulation comes from a point in the grid. Specifically, we use the range of the initial value to determine the sampling interval, divide the sampling interval into grids, and select points in the grid one by one to perform $N$ simulations, and the step size of each simulation is $L$.

% See numerical solution is presented above for three different trajectories. ${ }^{[4]}$

\subsection{Neural network training}
We can convert the real moment information into sigma points through unscented transformation, use neural ODE to fit the propagation function of sigma points, and get the predicted moment information from the predicted sigma points. Then calculate the training error and loss function to train the neural network.
% It is an effective way to transform moment information into sigma points through unscented transformation, and then estimate moment information from sigma points.We can fit the propagation function of the sigma point through the neural ODE, and then predict the sigma point at the next moment by calculating the numerical solution of the neural ODE, and then combine the supervised data to realize the training process of the neural network.

Assume that the mean $m \in \mathbb{R}^n$ and covariance $P \in \mathbb{R}^{n \times n}$ of the Monte Carlo simulation results at time $t$ are known. We can calculate the sigma point by Eq.\eqref{eq8}.
\begin{eqnarray}    \label{eq8}
&&z^{(0)}=m, \nonumber \\
&&z^{(i)}=m+[\sqrt{n+\lambda P}]_i, i=1, \ldots, n, \\
&&z^{(i)}=m-[\sqrt{n+\lambda P}]_{i-1}, i=n+1, \ldots, 2n, \nonumber 
\end{eqnarray}

Where $n$ is the state space dimension, $\lambda$ is the scaling factor, and $[A]_i$ is the column $i$ of the matrix $A$. We calculate the weight of each sigma point with respect to different moments by Eq.\eqref{eq9}.
\begin{eqnarray}    \label{eq9}
&&W_0^{(c)}=\frac{\lambda}{(n+\lambda)-\left(1-\alpha^2+\beta\right)} \nonumber \\
&&W_0^{(m)}=\frac{\lambda}{n+\lambda}  \\ 
&&W_i^{(c)}=W_i^{(m)}=\frac{1}{2(n+\lambda)}, i=1, \ldots, 2n \nonumber
\end{eqnarray}

Typically, we should set $\alpha$ to be small (e.g., $10^{-3}$), $\beta=2$ and $\lambda=1$ based on observations from Ref.\cite{julier1997new}. Through Eq.\eqref{eq8} and Eq.\eqref{eq9} we can get $W(t_k)^{m}$, $W(t_k)^{c}$ and $Z(t_k)$ at time $t=k$. We use the NeuralODE as a nonlinear map $F(z;\theta )$ of sigma points, satisfying the dynamical equation $\dot{z} =F(z;\theta )$. When $Z(t_k)$ is known, we can get $Z(t_{k+1})$ by integrating the NeuralODE. Then the mean and covariance at time $t=k+1$ can be predicted using Eq.\eqref{eq11} and Eq.\eqref{eq12}.
% We use the neural network to establish a nonlinear map $F$ to approximate the propagation function of the sigma point. 
% When the mean $m\left(t_k\right)$ and covariance $P\left(t_k\right)$ at time $t_k$ are known, the mean value $m\left(t_{k+1}\right)$ and covariance $P\left(t_{k+1}\right)$ can be obtained by integrating the differential Eq.\eqref{eq10} and Eq.\eqref{eq11} \cite{sarkka2007unscented}.
% \begin{eqnarray} 
%     \frac{d m(t)}{d t} &=&F(Z(t) ; \theta) w_m \label{eq10} \\
%     \frac{d P(t)}{d t} &=&Z(t) F^T(z(t) ; \theta)+F(z(t) ; \theta) W Z^T(t)+D Q_c(t) D^T  \label{eq11} 
% \end{eqnarray}
% , satisfying the formula $Z_{t+1}=F(Z_{t})$. Then the mean and covariance at time $t+1$ can be predicted using Eq.\eqref{eq11} and Eq.\eqref{eq12}.
\begin{eqnarray}
    &&\hat{\mu}(t_{k+1})=\sum_{i=0}^{2 n} w_i^{(m)} z(t_{k+1})^{(i)} \label{eq11}\\
    &&\widehat{\Sigma}(t_{k+1})=\sum_{i=0}^{2 n} w_i^{(c)}\left(z^{(i)}(t_{k+1})-\hat{\mu}(t_{k+1})\right)\left(z^{(i)}(t_{k+1})-\hat{\mu}(t_{k+1})\right)^T \label{eq12}
\end{eqnarray}
% \subsection{Train neural network}
% According to Simo Särkkä's research,

Because we assume that the disturbance of the system is Gaussian white noise with a mean of zero, two independent neural networks $F(z;\theta_1)$ and $F(z;\theta_2)$ can be used to predict the mean and covariance respectively. In particular, we sequentially solve the following two smaller optimization problems: 
\begin{eqnarray}
    &&\theta _{1}^{\star } =\underset{\theta_1 }{argmin} \sum_{k=0}^{L-1} \left \| \mu (t_{k+1}) -\hat{\mu }(t_{k+1}|t_{k};\theta _{1} ) \right \| ^{2},  \label{eq13}\\
    &&\theta _{2}^{\star } =\underset{\theta_2 }{argmin} \sum_{k=0}^{L-1} \left \| \Sigma (t_{k+1}) -\widehat{\Sigma }(t_{k+1}|t_{k};\theta _{1}^{\star },\theta _{2} ) \right \| ^{2}.  \label{eq12}
\end{eqnarray}

% Estimated values of $m\left(t_{k+1}\right)$ and $P\left(t_{k+1}\right)$ can be obtained by simple transformation:
% \begin{equation}
%     \hat{\mu}_x\left(t_{k+1} \mid t_k ; \theta\right)=\left[\begin{array}{ll}
% I & 0
% \end{array}\right] m\left(t_{k+1} ; \theta\right), \hat{\Sigma}_x\left(t_{k+1} \mid t_k ; \theta\right)=\left[\begin{array}{cc}
% I & 0 \\
% 0 & 0
% \end{array}\right] P\left(t_{k+1} ; \theta\right)
% \end{equation}

% The loss function is calculated by the following formula to train the nonlinear transformation $F$ parameterized by the neural network:
% \begin{equation}
%     L(\theta)=\sum_{k=1}^N\left\|\hat{\mu}_x\left(t_{k+1} \mid t_k ; \theta\right)-\mu_x\left(t_{k+1}\right)\right\|^2+\left\|\hat{\sum}_x\left(t_{k+1} \mid t_k ; \theta\right)-\sum_x\left(t_{k+1}\right)\right\|^2
% \end{equation}
% \subsection{Error Calculation}
We compute the probability density function(PDF) of the residuals using the  formula $p(\bm{x})=\varepsilon (\bm {x})/A$, where $\bm {x}$ belongs to the initial state space $\Omega \subset R^n$ and $A=\int _{\Omega } \varepsilon (\bm {x})dx$ is a normalizing constant. We use root mean square error (RMSE) and relative root mean square error (RRMSE) to evaluate the training results of the network, which can show different scales of error.

\subsection{Residual based sampling}

We assume that $p(x)$ is a second-order differentiable continuous function (if not we can use a polynomial function to fit $p(x)$). With the definition of maximum points and the continuity of $p(x)$, the set of maximum points $\omega = \{ x|p'(x)=0\ and\ p''(x)<0 \}$ is easy to find. The residual-based multi-peak sampling algorithm (RBMS) samples in the neighborhood $U(x_{i},r)$ of maximum point $x_{i}\in \omega$. The specific sampling process is shown in Algorithm \ref{AL1}.

\begin{algorithm}[!ht]
    \renewcommand{\algorithmicrequire}{\textbf{Input:}}
	\renewcommand{\algorithmicensure}{\textbf{Output:}}
	\caption{RBMS-I Algorithm with Sampling Hyperparameters}
    \label{power}
    \begin{algorithmic}[1] % 控制是否有序号 
    \label{AL1}
        \REQUIRE PDF of the residual $p(x)$; Hyperparameters $m,n,r$; Sample set $\Gamma _{k}$;% input 的内容
	    \ENSURE Sample set of the $(k+1)$th epoch $\Gamma _{k+1}$; % output 的内容
        \STATE $\omega \gets \{ x|p'(x)=0 \wedge  p''(x)<0 \wedge x\in \Omega \}$;
        % \STATE All the samples in the sample pool are sent to the neural network for training;
        % \STATE Calculate error $E_{mn}(\tau ), RE_{mn}(\tau )$ and residual $RMSE(\tau ), RRMSE(\tau )$;
        
        \STATE $\omega_n \gets \underset{\omega'\subset \omega,|\omega'|=n}{argmax}\{\sum p(x)|x\in \omega'\}  $
        
         % for loop
        \FOR {$x_{i}\in \omega_n$}
        
            \STATE $A\gets U(x_{i},r )\cap \Omega \setminus \Gamma _{k}$;

            \STATE $B\gets \underset{A'\subset A,|A'|=m}{argmax}\{\sum p(x)|x\in A'\}  $

            \STATE $\Gamma _{k}\gets \Gamma _{k}\cup B$
            
        \ENDFOR
        \STATE $\Gamma _{k+1}\gets \Gamma _{k}$
        
        % % if ... else
        % \IF {$RMSE(\tau ) \ge  RMSE(S)$ \textbf{or} $RRMSE(\tau ) \ge  RRMSE(S)$}
        %     \STATE Select $M$ non-adjacent maxima from the residual matrix $E_{mn}(\tau )$ or $RE_{mn}(\tau )$;
        %     \STATE For each maximum value point, select $N$ points closest to it;
        %     \STATE A total of $MN$ sample point sample set $S_{MN}$, satisfying $S_{MN} \cap \tau =\phi$;
        %     \STATE $\tau \gets \tau \cup S_{MN}$;
        % \ELSE
        %     \BREAK;
        % \ENDIF
        
        % while 
        % \WHILE {$RMSE(\tau ) \ge  RMSE(S)$ \textbf{or} $RRMSE(\tau ) \ge  RRMSE(S)$}
        %     \STATE Select $M$ non-adjacent maxima from the residual matrix $E_{mn}(\tau )$ or $RE_{mn}(\tau )$;
        %     \STATE For each maximum value point, select $N$ points closest to it;
        %     \STATE A total of $MN$ sample point sample set $S_{MN}$, satisfying $S_{MN} \cap \tau =\phi$;
        %     \STATE $\tau \gets \tau \cup S_{MN}$;

        %     \STATE All the samples in the sample pool are sent to the neural network for training;

        %     \STATE Calculate error $E_{mn}(\tau ), RE_{mn}(\tau )$ and residual $RMSE(\tau ), RRMSE(\tau )$;
            
        % \ENDWHILE
        
        % % do ... while
        % \REPEAT
        %     \STATE a--
        % \UNTIL {a < 0}
        
        % \STATE \textbf{return} $\mathbf{Q}_{Iter}$.
    \end{algorithmic}
\end{algorithm}

The training set of the first epoch of the neural network training process comes from several sample points randomly selected in the sample space $\Omega$. Then we repeatedly execute the RBMS algorithm to generate the data required for each epoch in the neural network training process. Such repeated operations will not stop until the error (RMSE, RRMSE, etc.) of the neural network training result is less than a preset threshold.

Compared with the method of collecting $k$ points with the largest residual error, the RBMS algorithm takes into account the influence of multiple maximum values. The RBMS algorithm can automatically detect the region where the maximum value exists, and divide the subspace that needs to be sampled in the entire sampling space. Compared with the fixed subspace segmentation method, the RBMS algorithm is more adaptive.

% We first used Algorithm \ref{AL1} for a case study, and soon we found two problems: 1) The sampling hyperparameters M and N are difficult to determine. 2) The number of peaks in the residual matrix is different in different training generations. Therefore, we decided to abandon the two hyperparameters of M and N, and sample a point at each peak. The process is shown in Algorithm \ref{AL2}.
We found that the number of maximum points of the function $p(x)$ is different in different training epochs when using the RBMS-I algorithm for the case study. Determining the two hyperparameters $m$ and $n$ is difficult. Therefore, we decided to abandon the two hyperparameters $m$ and $n$, set $r$ to a number greater than 1 (for example $r=$5), and sample one point at each maximum value. We name the algorithm without hyperparameters as the RBMS-II algorithm, and the process is shown in Algorithm \ref{AL2}.

\begin{algorithm}[!ht]
    \renewcommand{\algorithmicrequire}{\textbf{Input:}}
	\renewcommand{\algorithmicensure}{\textbf{Output:}}
	\caption{RBMS-II Algorithm without Sampling Hyperparameters}
    \label{power}
    \begin{algorithmic}[1] % 控制是否有序号 
    \label{AL2}
        \REQUIRE PDF of the residual $p(x)$; Sample set $\Gamma _{k}$;% input 的内容
	    \ENSURE Sample set of the $(k+1)$th epoch $\Gamma _{k+1}$; % output 的内容
        \STATE $\omega \gets \{ x|p'(x)=0 \wedge  p''(x)<0 \wedge x\in \Omega \}$;
        % \STATE All the samples in the sample pool are sent to the neural network for training;
        % \STATE Calculate error $E_{mn}(\tau ), RE_{mn}(\tau )$ and residual $RMSE(\tau ), RRMSE(\tau )$;
        
        % \STATE $\omega_n \gets \underset{\omega'\subset \omega,|\omega'|=n}{argmax}\{\sum p(x)|x\in \omega'\}  $
        
         % for loop
        \FOR {$x_{i}\in \omega_n$}
        
            \STATE $A\gets U(x_{i},5 )\cap \Omega \setminus \Gamma _{k}$;

            \STATE $B\gets \underset{x}{argmax}\{ p(x)|x\in A\}  $

            \STATE $\Gamma _{k}\gets \Gamma _{k}\cup B$
            
        \ENDFOR
        \STATE $\Gamma _{k+1}\gets \Gamma _{k}$
        
        % % if ... else
        % \IF {$RMSE(\tau ) \ge  RMSE(S)$ \textbf{or} $RRMSE(\tau ) \ge  RRMSE(S)$}
        %     \STATE Select $M$ non-adjacent maxima from the residual matrix $E_{mn}(\tau )$ or $RE_{mn}(\tau )$;
        %     \STATE For each maximum value point, select $N$ points closest to it;
        %     \STATE A total of $MN$ sample point sample set $S_{MN}$, satisfying $S_{MN} \cap \tau =\phi$;
        %     \STATE $\tau \gets \tau \cup S_{MN}$;
        % \ELSE
        %     \BREAK;
        % \ENDIF
        
        % while 
        % \WHILE {$RMSE(\tau ) \ge  RMSE(S)$ \textbf{or} $RRMSE(\tau ) \ge  RRMSE(S)$}
        %     \STATE Select $M$ non-adjacent maxima from the residual matrix $E_{mn}(\tau )$ or $RE_{mn}(\tau )$;
        %     \STATE For each maximum value point, select $N$ points closest to it;
        %     \STATE A total of $MN$ sample point sample set $S_{MN}$, satisfying $S_{MN} \cap \tau =\phi$;
        %     \STATE $\tau \gets \tau \cup S_{MN}$;

        %     \STATE All the samples in the sample pool are sent to the neural network for training;

        %     \STATE Calculate error $E_{mn}(\tau ), RE_{mn}(\tau )$ and residual $RMSE(\tau ), RRMSE(\tau )$;
            
        % \ENDWHILE
        
        % % do ... while
        % \REPEAT
        %     \STATE a--
        % \UNTIL {a < 0}
        
        % \STATE \textbf{return} $\mathbf{Q}_{Iter}$.
    \end{algorithmic}
\end{algorithm}

The steps of RBMS-II algorithm are simpler than RBMS-I. More critically, the RBMS-II algorithm is a fully adaptive sampling algorithm that does not depend on any hyperparameters. This means that the RBMS-II algorithm can be easily applied to different tasks without spending effort on parameter tuning.

\begin{figure}[H]
    \centering
    \includegraphics[width=16.5cm]{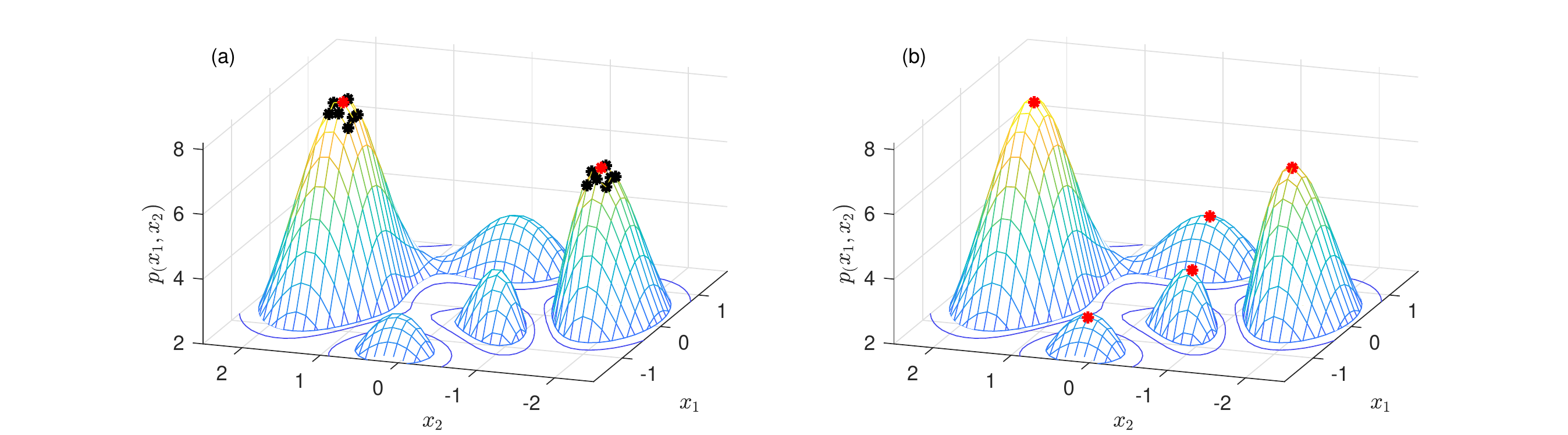}
    \caption{Display of RBMS algorithm sampling position.}
    \label{fig:0}
\end{figure}

Fig.\ref{fig:0}(a) shows the sampling situation of the RBMS-I algorithm when $m=7, n=2, r=5$. The sampling of the RBMS-II algorithm is shown in Fig.\ref{fig:0}(b). The mark $\ast$ indicates the sampling position, the red mark represents the maximum value point, and the black mark represents the sampling point near the maximum value point.

% \subsection{Analyzing Dynamic Behavior}

% Stochastic dynamical systems have rich dynamic behaviors, such as P bifurcation, D bifurcation, first crossing probability, maximum possible path, etc., which are often studied.

% 由 1.3.2节随机吸引域的定义可知, 为了得到不同草食动物密度 $H$ 和高斯噪声强度 D 下高植被量随机吸引域 我们首先要计算植被从确定性高植被域 $D$ 经过边界 $\Gamma_1$ 到 达确定性低植被域 $D^c$ 的逃逸概率 $P_E(x)$ 来确定 $D_I$, 然后计算从 $D_I^c$ 经过边界 $\Gamma_{1 I}$ 到 达 $D_I$ 的逃逸概率 $P_R(x)$ 来确定 $D_{I I^{\circ}}^c$ 根据文献 ${ }^{[109]}$, 对于方程 $(2-5), P_E(x)$ 和 $P_R(x)$ 是下列 Dirichlet 边界值问题的解
% $$
% A_1 P_E(x)=0
% $$
% 边界条件为
% $$
% \underline{\left.P_E\right|_{\Gamma_1}=1,},\left.P_E\right|_{\Gamma_2}=0
% $$
% 以及
% $$
% A_1 P_R(x)=0
% $$
% 边界条件为
% $$
% \left.P_R\right|_{\Gamma_{1 I}}=1,\left.\quad P_R\right|_{\Gamma_{2 I}}=0,
% $$
% 这里 $A_1$ 是方程 $(2-5)$ 的生成算子
% $$
% A_1 P=m(x) \frac{d P}{d x}+\frac{1}{2} \sigma^2(x) \frac{d^2 P}{d^2 x} \cdot \quad \begin{aligned}
% & \text { MFET } \\
% & =-1
% \end{aligned}
% $$

\section{Result}

\subsection{A Model of Vegetation Biomass Change in Grazing Systems}
In general, changes in vegetation biomass depend on the net growth rate of vegetation and the rate of vegetation consumption. The net growth rate of vegetation can be expressed by the well-known Logistic function, as shown in Eq.\eqref{18}. Vegetation consumption rate can be represented by a saturation function in the form of Eq.\eqref{19}.

\begin{equation}
    G(x)=kx(A-x) \label{18}
\end{equation}

\begin{equation}
    R(x)=\beta \frac{cx^{2} }{x^{2}+x_{0}^{2}}  \label{19}
\end{equation}

As shown in Eq.\eqref{20}, a classical vegetation biomass dynamics model can be obtained by combining the above two equations  \cite{noy1975stability}.

\begin{equation}
    \dot{x}=G(x)-R(x)=kx(A-x) - \beta \frac{cx^{2} }{x^{2}+x_{0}^{2}}  \label{20}
\end{equation}

where $k$ is the growth coefficient, $A$ is the environmental holding capacity, $\beta$ is the herbivore density, $c$ is the consumption coefficient, and $x_{0}$ corresponds to the vegetation biomass when the grazing rate reaches half the maximum.

This system has two steady states\cite{zhang2019noise}. There are inevitable random factors in the real ecological environment, so we consider adding Gaussian white noise to the model, as shown in Eq.\eqref{21}.
\begin{equation}
    \dot{x}=kx(A-x) - \beta \frac{cx^{2} }{x^{2}+x_{0}^{2}} +x\xi (t) \label{21}
\end{equation}

Because noise is considered in the system, the solution of the system may jump between two attractive domains, and this jump may cause the system to tend to a different equilibrium state, which will have a huge impact on people's production. Therefore, it is valuable to fit the dynamic model and analyze its first crossing probability.

As mentioned earlier, we first convert the equation of the form $\dot{X}=F(X,t;\xi)$ into a stochastic differential equation of the Ito form. By adding the Wong-Zakai correction term\cite{eugene1965relation}, the Eq.\eqref{21} can be transformed into an It$ \hat{\rm o}$ stochastic differential equation, as shown in Eq.\eqref{22}.

\begin{equation}
    dX=m(X)dt+\sigma (X)dB(t) \label{22}
\end{equation}

Where $B(t)$ is a unit Wiener process, and the drift and diffusion terms are defined as shown in Eq.\eqref{23}. Suppose Gaussian noise intensity is $D$.
% \begin{eqnarray}
%     m(X)&=&kx(A-x) - \beta \frac{cx^{2} }{x^{2}+x_{0}^{2}} + DX \label{23} \\
%     \sigma (X)&=&\sqrt{2D} X \nonumber
% \end{eqnarray}
\begin{equation} \label{23}
\begin{aligned}
&  m(X)=kx(A-x) - \beta \frac{cx^{2} }{x^{2}+x_{0}^{2}} + DX \\
& \sigma (X)=\sqrt{2D} X
\end{aligned}
\end{equation}

We numerically solve Eq.\eqref{22} using the Euler-Milstein method in the form of Eq.\eqref{5}. The solution step is 0.01, and the number of steps is 25. In order to count the stable moments, we simulated the above solution process 10,000 times using the Monte Carlo simulation method \cite{o2022stochastic}.

\begin{figure}
    \centering
    \includegraphics[width=7cm]{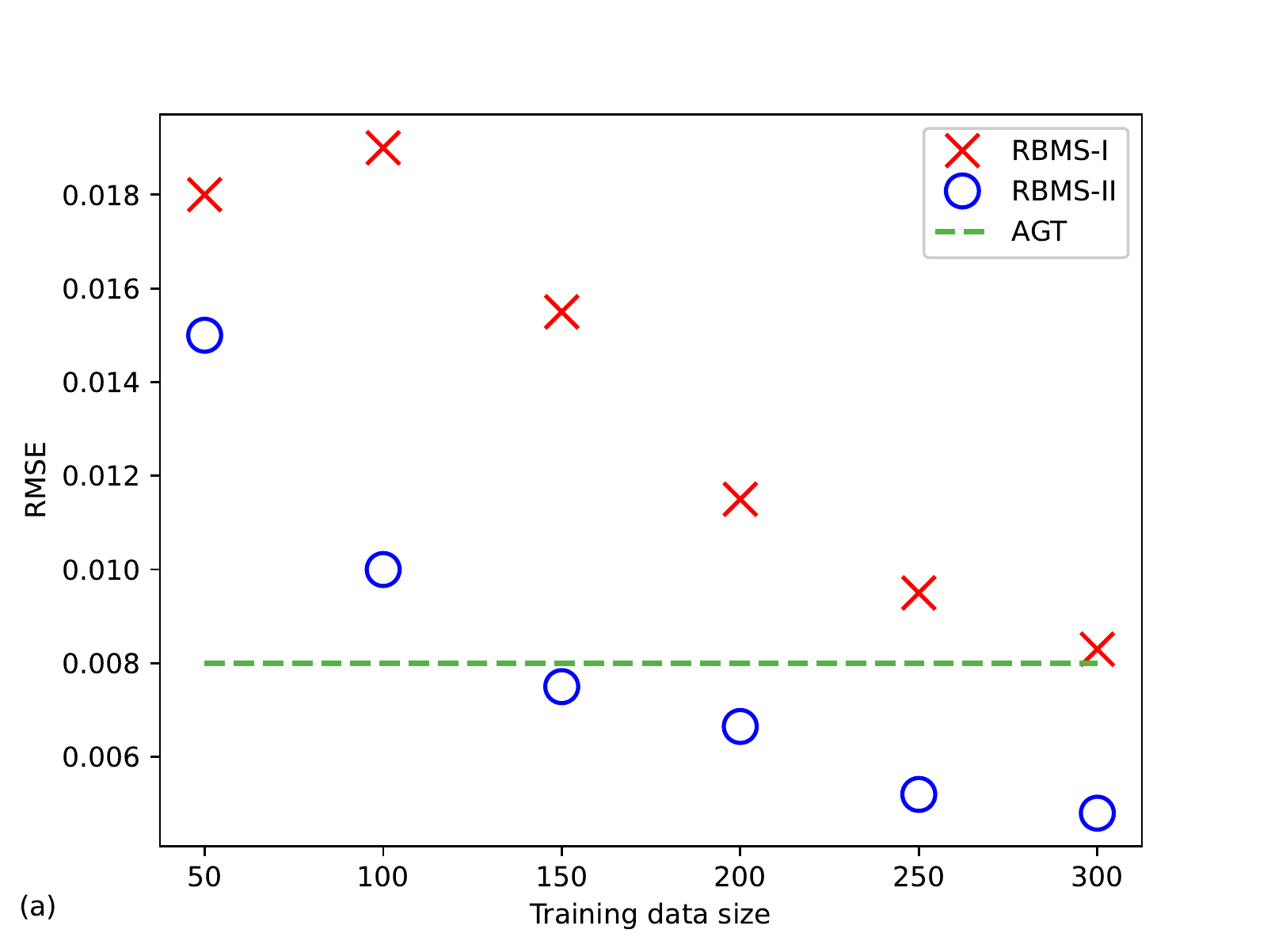}
    \includegraphics[width=7cm]{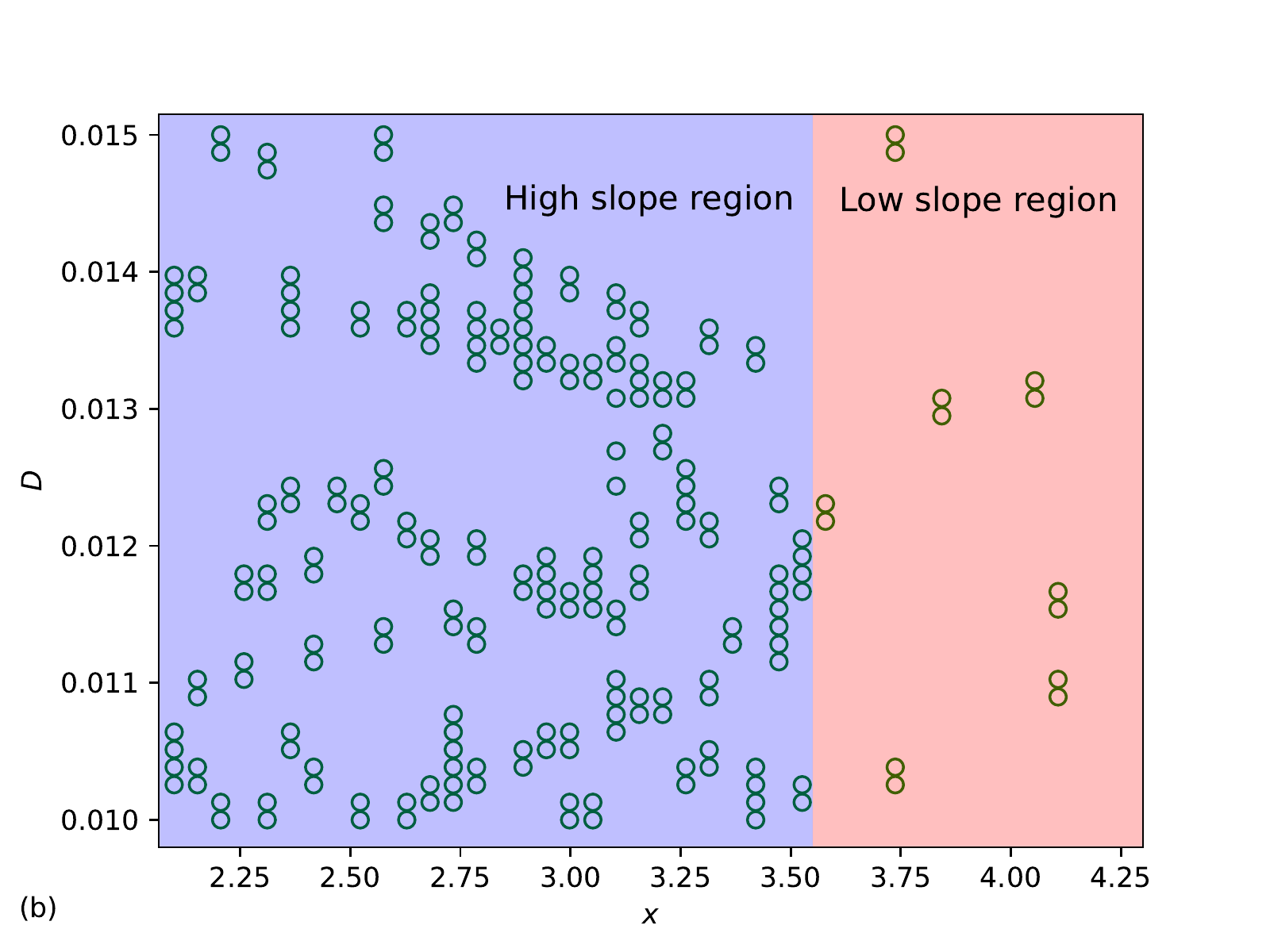}
    \caption{Analysis of error and sampling distribution after fitting grazing system with RBMS algorithm.}
    \label{fig:1}
\end{figure}

In the Fig.\ref{fig:1}(a), we show the training effects of the three point-taking strategies. The three point-taking strategies are: all points in the grid (AGF, 1600), RBMS-I algorithm and RBMS-II algorithm. We use the training error (RMSE) of all sample points as the threshold. The number of samples selected by the RBMS algorithm is less than 300 so that the training error reaches the threshold, which means that using all points for training is very inefficient. To achieve the same effect RBMS-II uses fewer sample points (about 150) than RBMS-I algorithm (about 300). RBMS-II algorithm is obviously better than RBMS-I algorithm.

In the Fig.\ref{fig:1}(b), we show the specific sampling positions in a fitting task. We first notice that the sampling is not uniform under the guidance of the RBMS algorithm. In the entire sampling space, there are densely sampled regions and sparsely sampled regions. We try to find the relationship between the density of sampling and the function being fitted. We noticed that more samples are often needed in the region where the slope of the fitted function is larger, and only a few samples are needed in the region where the slope of the fitted function tends to zero. This may mean that the function is easier to learn where the slope tends to zero.

\begin{figure}[H]
    \centering
    \includegraphics[width=7cm]{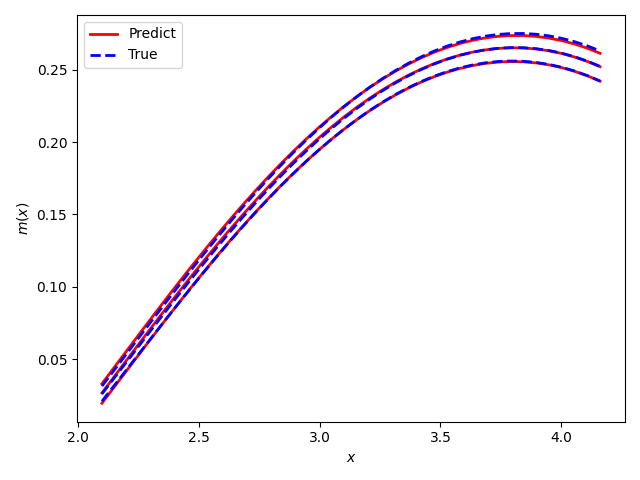}
    \includegraphics[width=7cm]{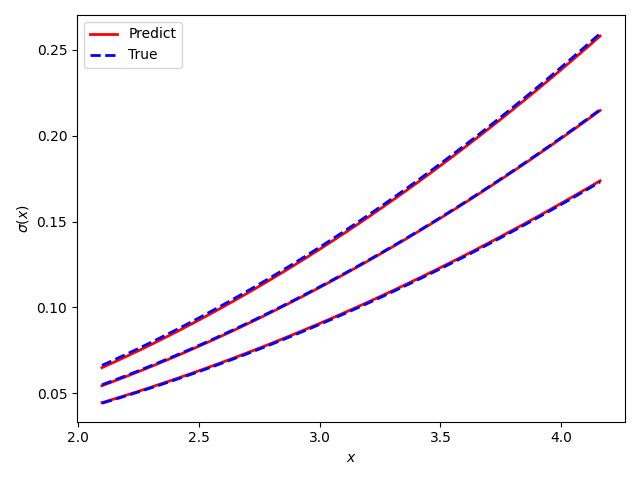}
    \caption{Learning drift and diffusion terms for grazing system using RBMS algorithms.}
    \label{fig:2}
\end{figure}

Although the main method used in the fitting task is not from us, the same fitting effect can be obtained with only a few samples after combining with RBMS. In Fig.\ref{fig:2}, we used three different noise intensities to make three curves, and we can see that the fitting effect of the neural network is very good.

To test the accuracy of the fitted models, we performed further kinetic analyze using surrogate models. Because the grazing system is a bistable model, we next analyze the escape probability $P_{E} (x)$ from a domain of high vegetation attraction to a domain of low vegetation attraction. For Eq.\eqref{22}, $P_{E} (x)$ is the solution of the following Dirichlet boundary value problem\cite{duan2015introduction}:
\begin{equation}
    A_1 P_E(x)=0, \left.P_E\right|_{\Gamma_1}=1, \left. P_E\right|_{\Gamma_2}=0. \label{24}
\end{equation}

Here $A_1$ is the generating operator for the Eq.\eqref{22}:
\begin{equation}
    A_1 P=m(x) \frac{d P}{d x}+\frac{1}{2} \sigma^2(x) \frac{d^2 P}{d^2 x} \cdot 
\end{equation}

By solving Eq.\eqref{24}, the analytical form of the escape probability can be obtained, as shown in Eq.\eqref{26}. Then use the difference method to calculate the numerical solution.

\begin{equation}
    P_E(x)=\int C_1 e^{\int \frac{-2 m(x)}{\sigma^2(x)} d x} d x+C_2 \label{26}
\end{equation}

\begin{figure}[H]
    \centering
    \includegraphics[width=5.3cm, trim=11 6 15 9,clip]{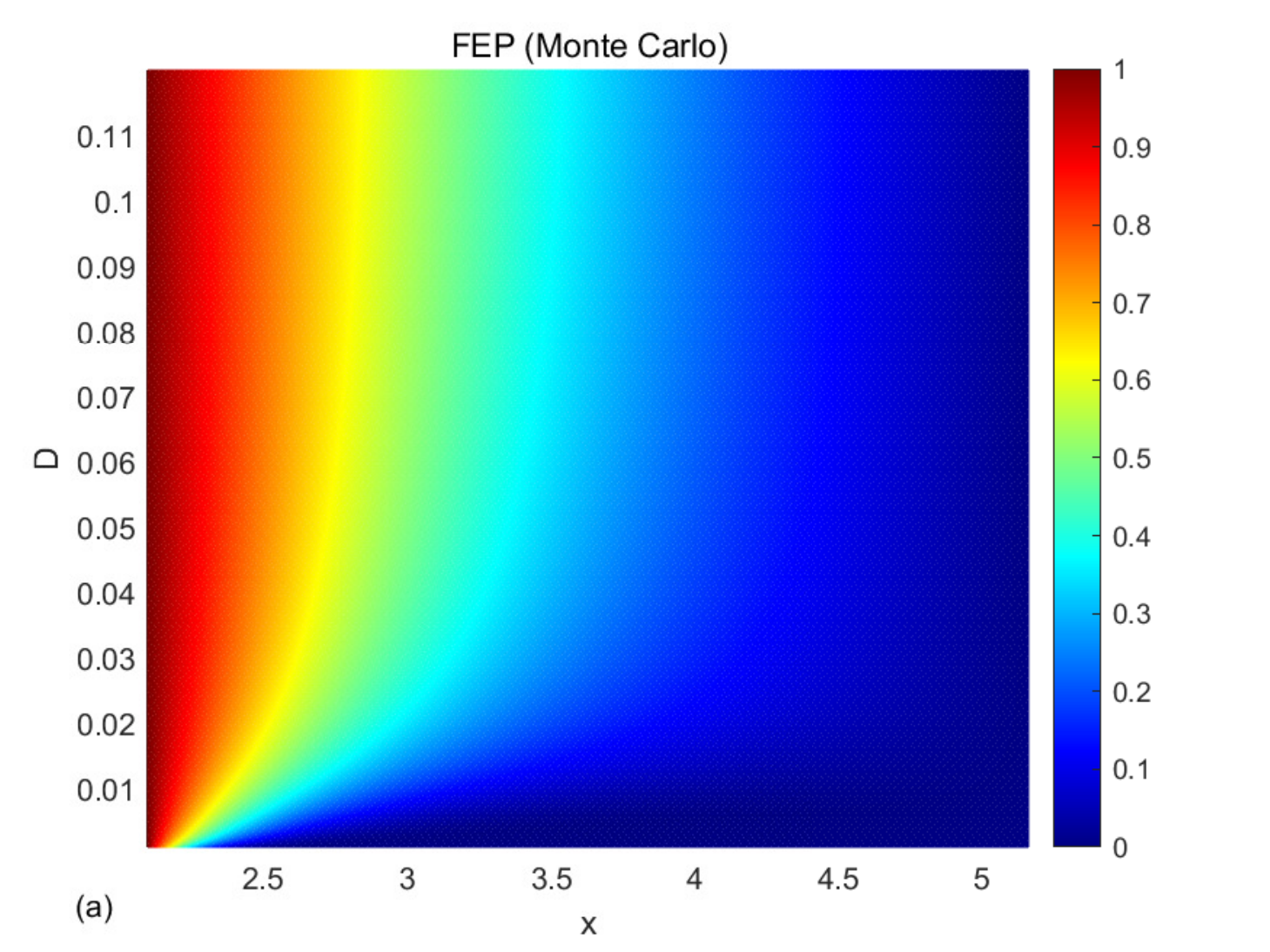}
    \includegraphics[width=5.3cm, trim=11 6 15 9,clip]{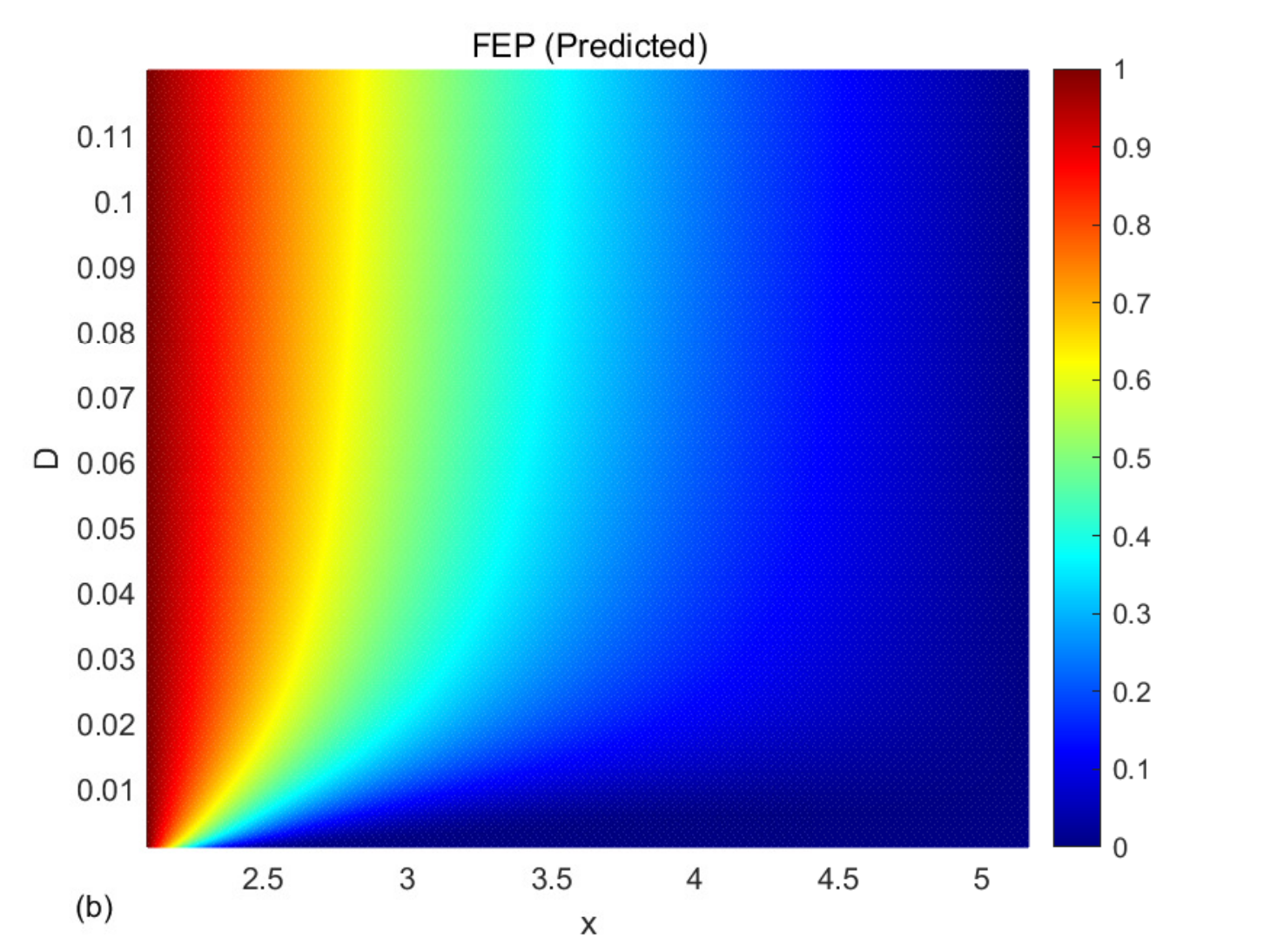}
    \includegraphics[width=5.3cm, trim=11 6 15 9,clip]{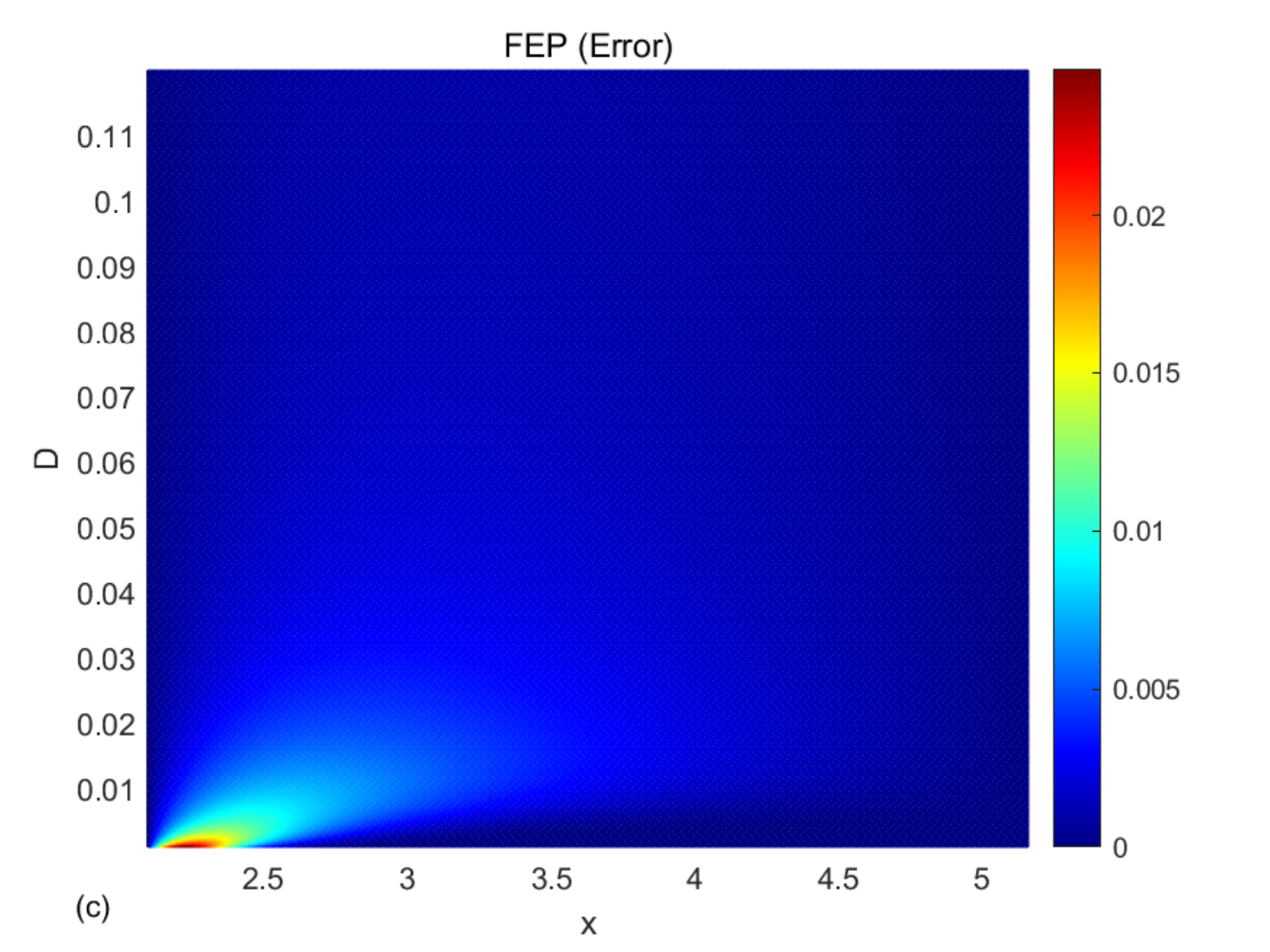}
    \caption{Top view of first escape probability density function based on analytical method and Monte Carlo simulation.}
    \label{fig:3}
\end{figure}

We can easily obtain the true escape probability density function by performing Monte Carlo simulation on Eq.\eqref{21}. $\theta _{1}(x)$ and $\theta _{2}(x)$ are the parameterized functions of the trained neural network. We set $m(x)=\theta _{1}(x), \sigma (x)=\theta _{2}(x)$, and then use Eq.\eqref{26} to solve the escape probability density function, which is the analytical method. The final results of these two methods are shown in Fig.\ref{fig:3}(a,b), and the two results are very close as shown in Fig.\ref{fig:3}(c). This shows that the fitting result after using the RBMS algorithm is reliable, and further dynamic behavior analysis can be carried out.

\subsection{Rayleigh-Van der pol impact vibration system}

In recent decades, in engineering mechanics, applied mathematics, applied physics and other fields, various collision systems have received extensive attention and research\cite{paget1937vibration,foale1994bifurcations,ibrahim2009vibro}. Introducing random noise excitations into the system may cause the system to exhibit stochastic bifurcation behavior. In a random impact vibration system, if the stable most probable steady-state trajectory changes from a point to a semicircle, or from a semicircle to a point, it means that the system has a P bifurcation phenomenon.

We consider a Rayleigh-Van der pol impact vibration system excited by multiplicative Gaussian white noise. The collision occurs at the position where the displacement is zero. The system equation is shown in Eq.\eqref{27}.

\begin{equation} \label{27}
\begin{aligned}
& \ddot{x}-\left(\alpha-\beta x^2\right) \dot{x}+\gamma \dot{x}^3+x=\dot{x} \xi_2(t), \quad x>0  \\
&  \dot{x}_{+}=-r \dot{x}_{-}, \quad x=0 
\end{aligned}
\end{equation}

% \begin{eqnarray}
%     \ddot{x}-\left(\alpha-\beta x^2\right) \dot{x}+\gamma \dot{x}^3+x&=&\dot{x} \xi_2(t), \quad x>0 \label{27} \\
%     \dot{x}_{+}&=&-r \dot{x}_{-}, \quad x=0 \nonumber
% \end{eqnarray}

where the parameters $\alpha, \beta$ and $\gamma$ are real numbers. $\xi_i(t)(i=1,2)$ is independent Gaussian white noise. By adding the Wong-Zakai correction term, the Eq.\eqref{27} can be transformed into an It$ \hat{\rm o}$ stochastic differential equation in the form shown in Eq.\eqref{28}.

\begin{equation} \label{28}
\begin{aligned}
& d x_1(t)=m_1\left(x_1, x_2\right) d t+\sigma_1\left(x_1, x_2\right) d B(t)  \\
& d x_2(t)=m_2\left(x_1, x_2\right) d t+\sigma_2\left(x_1, x_2\right) d B(t) 
\end{aligned}
\end{equation}

% \begin{eqnarray}
% d x_1(t)&=&m_1\left(x_1, x_2\right) d t+\sigma_1\left(x_1, x_2\right) d B(t) \label{28}
% \\
% d x_2(t)&=&m_2\left(x_1, x_2\right) d t+\sigma_2\left(x_1, x_2\right) d B(t) \nonumber
% \end{eqnarray}

Where $B(t)$ is a unit Wiener process, and the drift and diffusion terms are defined as shown in Eq.\eqref{29}. Suppose Gaussian noise intensity is $D$.
\begin{equation} \label{29}
\begin{aligned}
& m_1\left(x_1, x_2\right)=x_2  \\
& \sigma_1\left(x_1, x_2\right)=\sqrt{2D} x_{1}  \\
& m_2\left(x_1, x_2\right)=(\alpha-\beta x_{1}^2)x_{2}-\gamma x_{2}^3 - x_{1} + Dx_{2} \\
& \sigma_2\left(x_1, x_2\right)=\sqrt{2D} x_{2}
\end{aligned}
\end{equation}

% \begin{eqnarray}
%     m_1\left(x_1, x_2\right)&=&x_2 \nonumber \\
%     \sigma_1\left(x_1, x_2\right)&=&0  \label{29} \\
%     m_2\left(x_1, x_2\right)&=&(\alpha-\beta x_{1}^2)x_{2}-\gamma x_{2}^3 - x_{1} + Dx_{2}  \nonumber \\
%     \sigma_2\left(x_1, x_2\right)&=&\sqrt{2D} x_{2} \nonumber
% \end{eqnarray}

In the numerical experiment part, we use the Euler-Milstein method to solve Eq.\eqref{28}, set the step size to 25 steps, $dt=0.01$. In the sampling space of [1,3]$\times $[-3,-1], we uniformly sample 1600 (40$\times$40) initial values are used as all point training.

\begin{figure}[H]
    \centering
    \includegraphics[width=7cm]{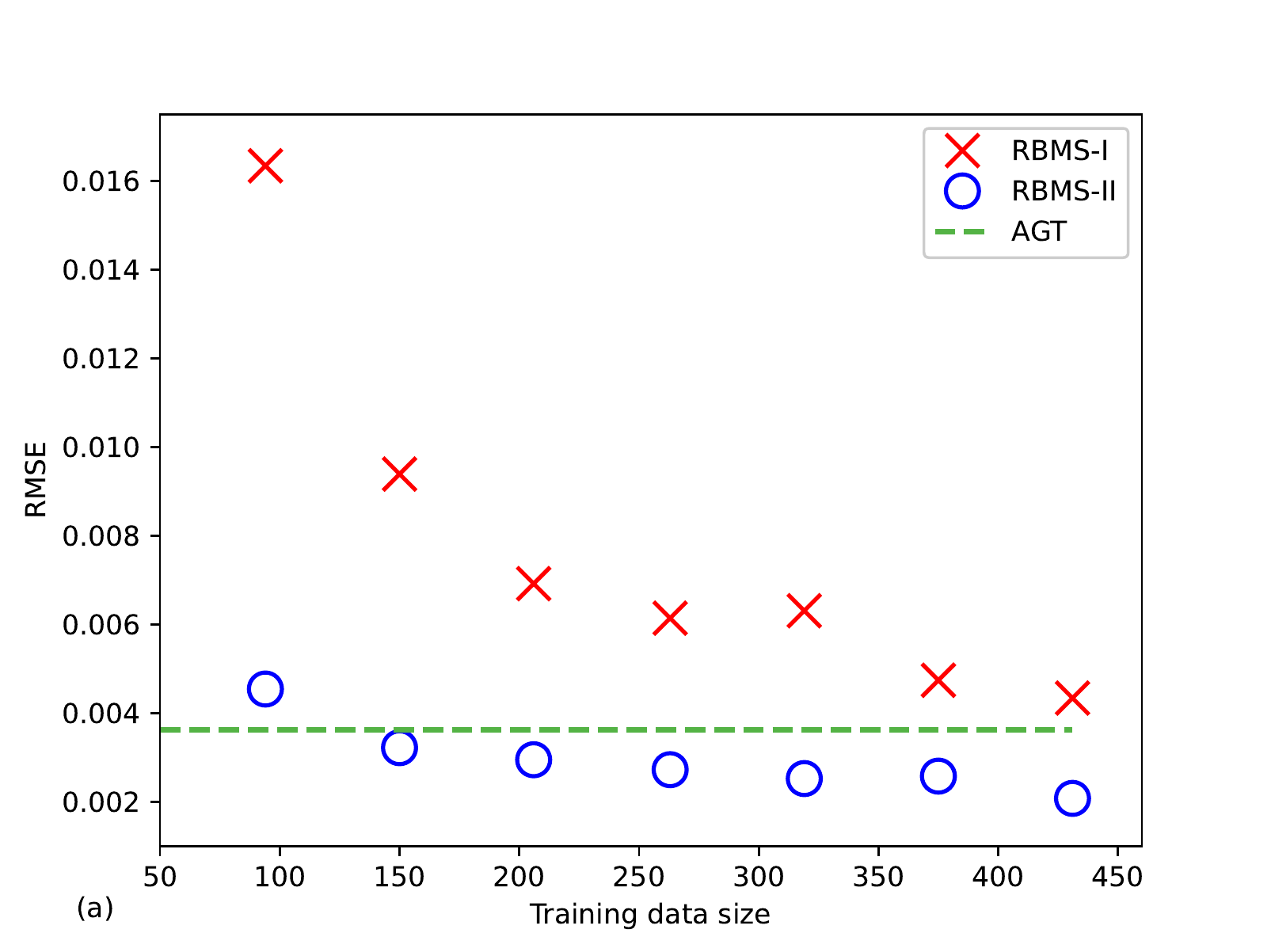}
    \includegraphics[width=7cm]{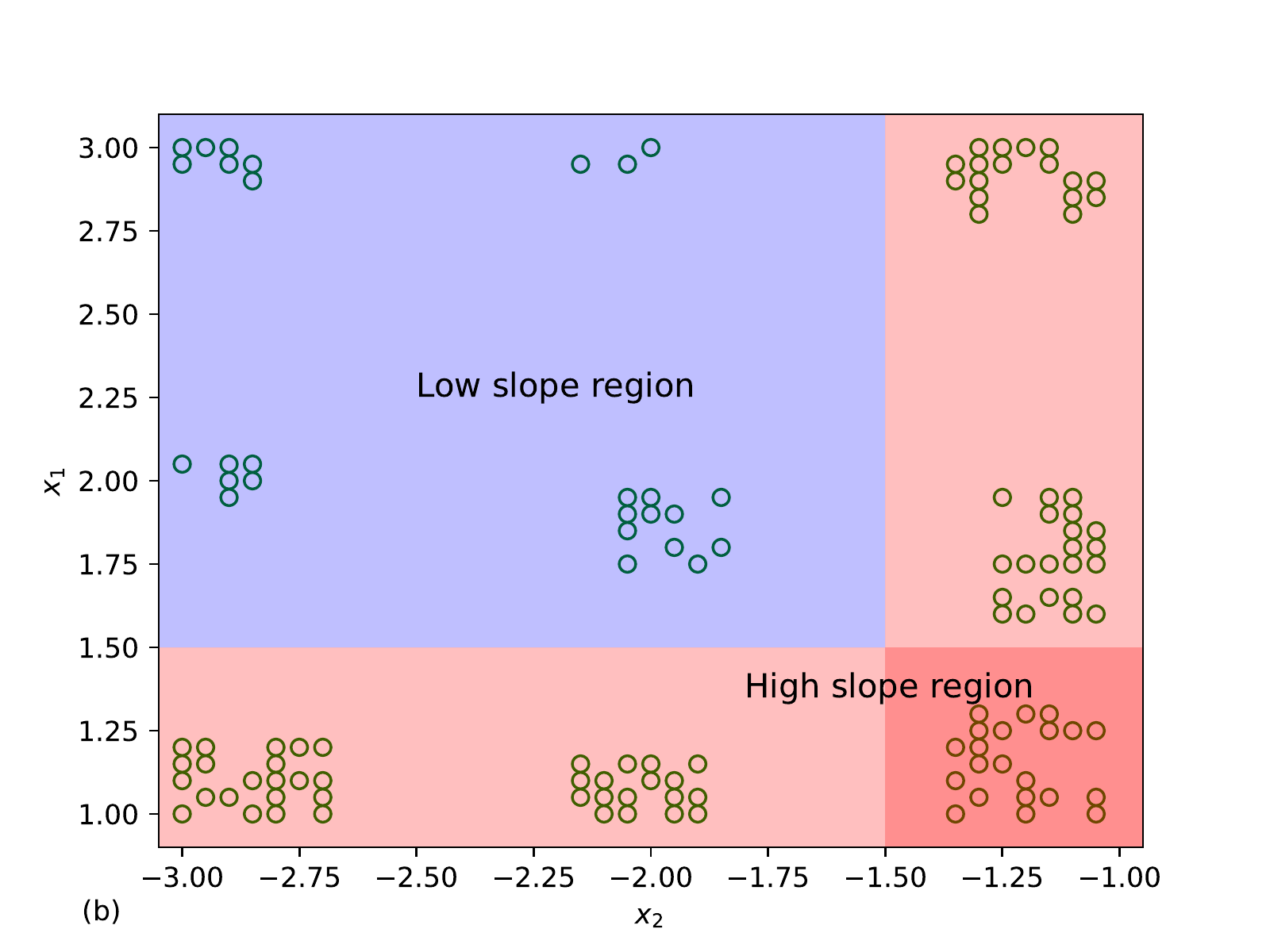}
    \caption{Analysis of error and sampling distribution after fitting Rayleigh-Van der pol system with RBMS algorithm.}
    \label{fig:4}
\end{figure}

In the Fig.\ref{fig:4}(a), we use the training results of all sample points as a benchmark to compare the error reduction of the two RBMS algorithms. Observing the image shows that the RBMS-II algorithm falls below the baseline error faster than the RBMS-I algorithm, although the RBMS-II algorithm can also reach the baseline error at about 450 sampling points. This shows that the RBMS-II algorithm is a more efficient sampling algorithm.

In the Fig.\ref{fig:4}(b), we draw the specific coordinates of the sampling points to visually observe the distribution of the sampling points. We have observed the results of many experiments and obtained two general rules: 1) The RBMS algorithm tends to sample at the boundary of the sampling space. 2) The RBMS algorithm is more inclined to sample in the area with a high slope of the fitted function. This is consistent with what was observed in the case of grazing systems. Therefore, we believe that the boundary and high slope areas are the key areas for sampling.

\begin{figure}[htb]
    \centering
    \includegraphics[width=4.5cm]{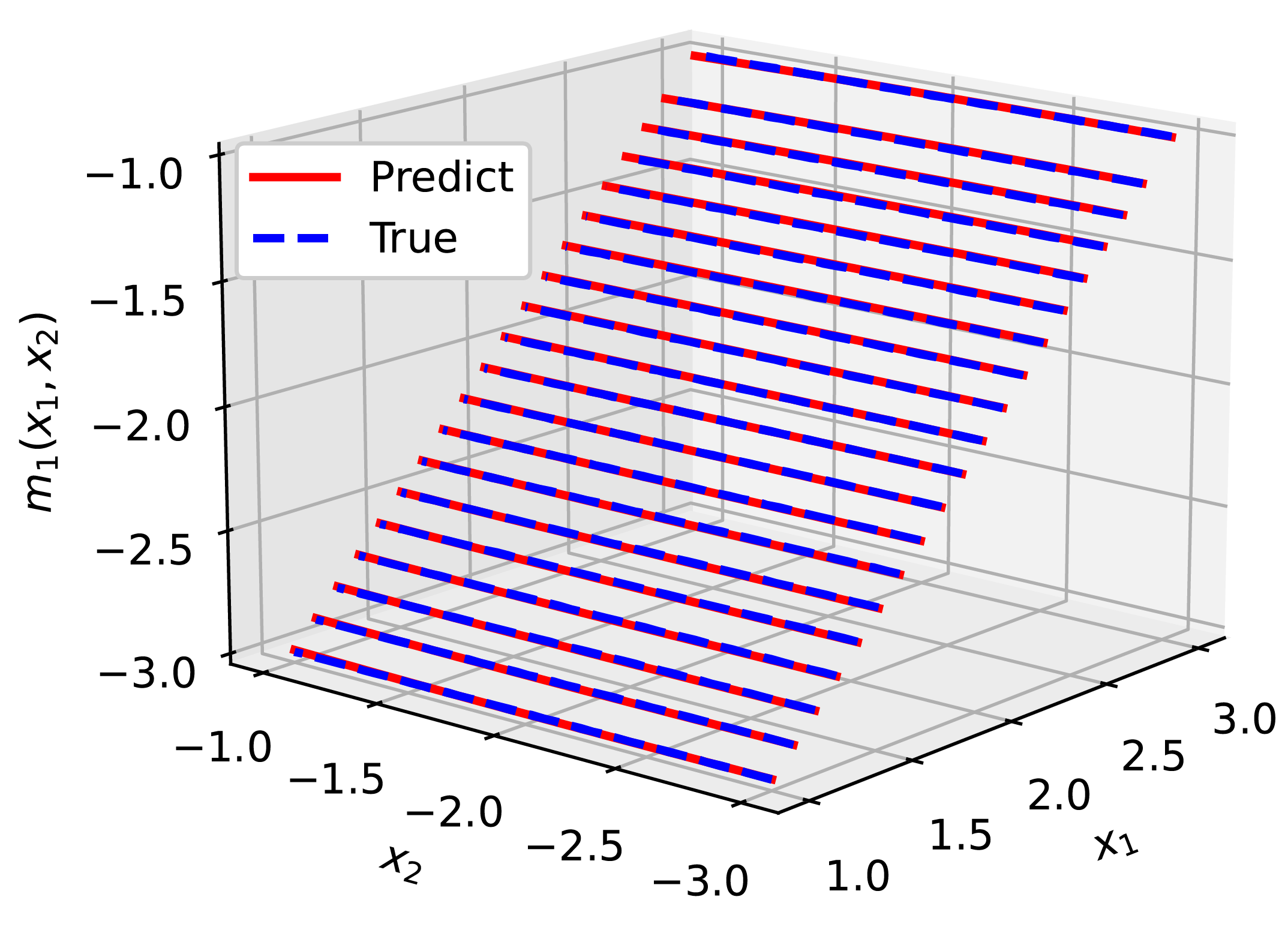}
    \includegraphics[width=4.5cm]{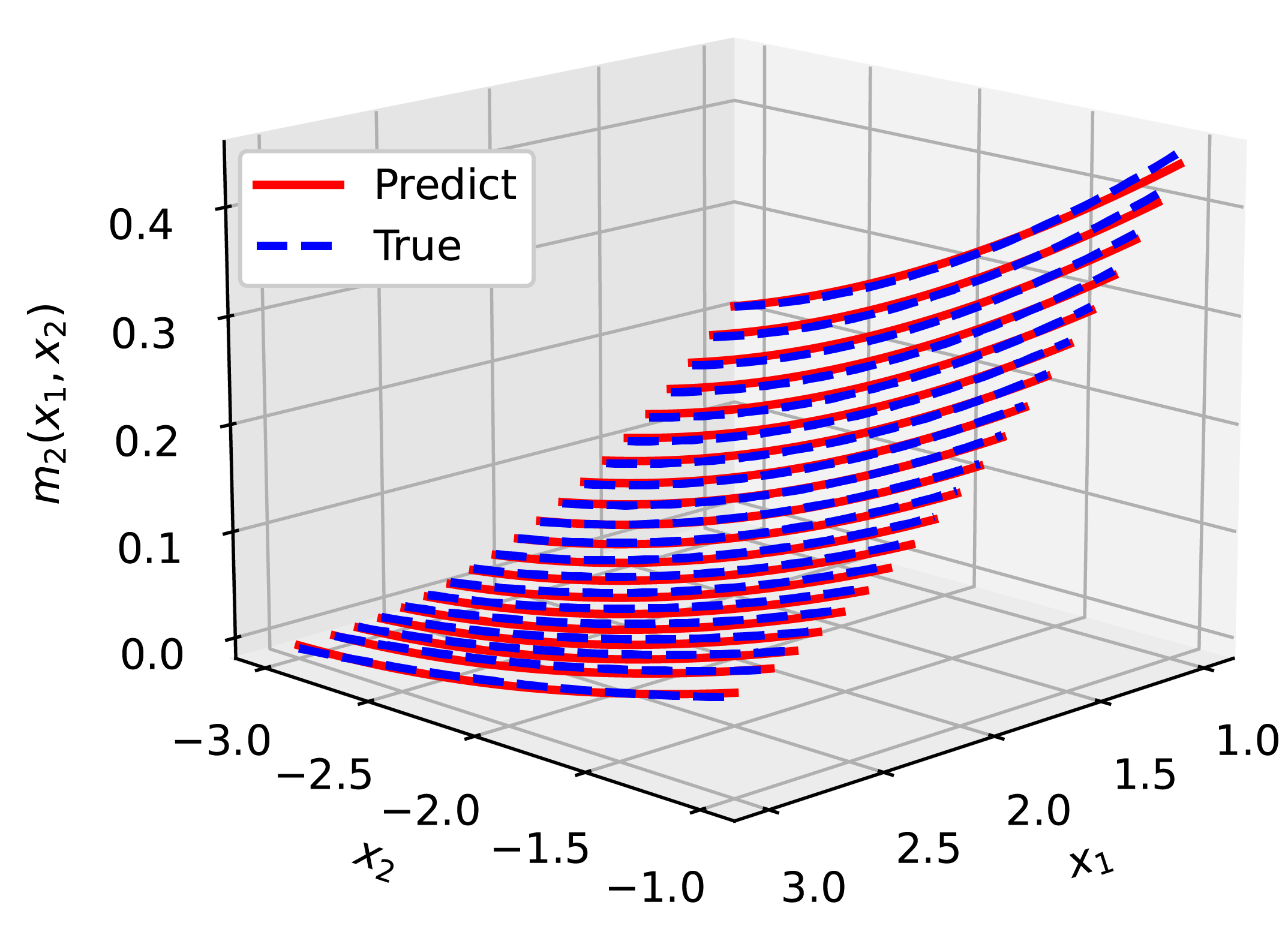} \\
    \includegraphics[width=4.5cm]{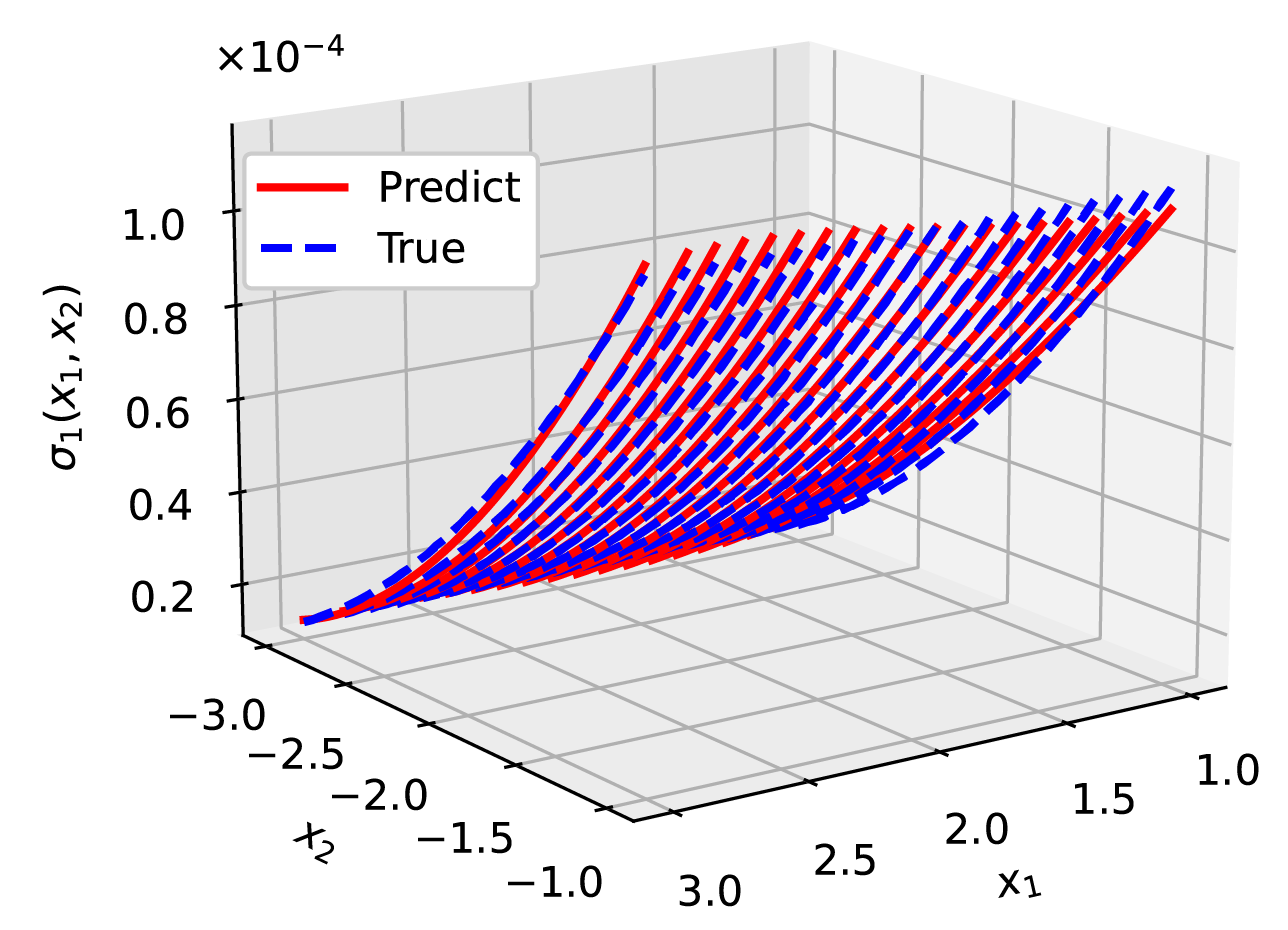}
    \includegraphics[width=4.5cm]{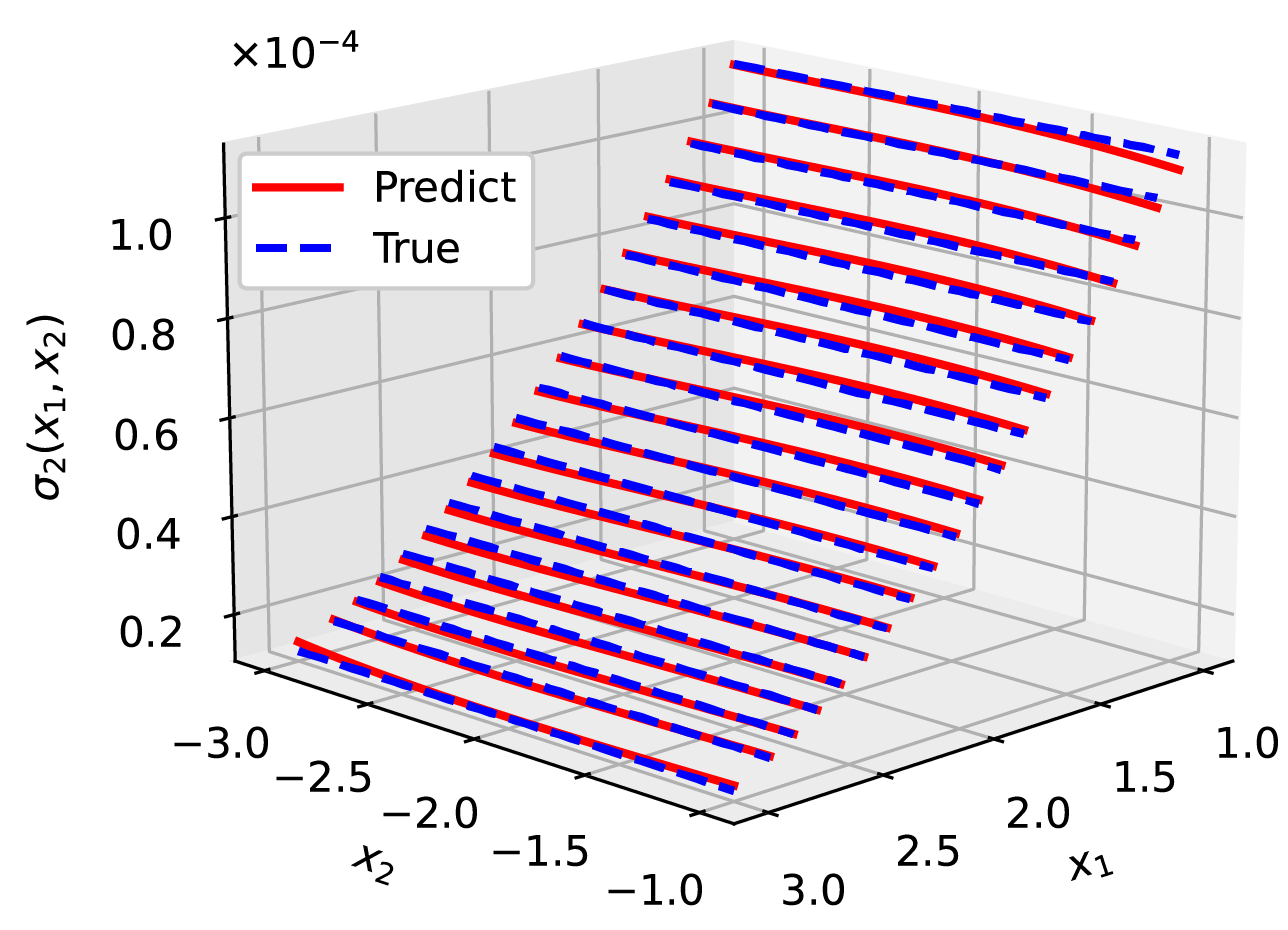}
    \caption{Learning drift and diffusion terms for Rayleigh-Van der pol system using RBMS algorithms.}
    \label{fig:5}
\end{figure}

In Fig.\ref{fig:5}, we compare the values of the drift and diffusion terms predicted by the neural network with the true values. Using the RBMS algorithm can achieve accurate fitting with fewer sampling points. Based on such results, we are confident to use neural network functions instead of drift and diffusion terms for further dynamical analysis.

\begin{figure}[htb]
    \centering
    \includegraphics[width=5.3cm]{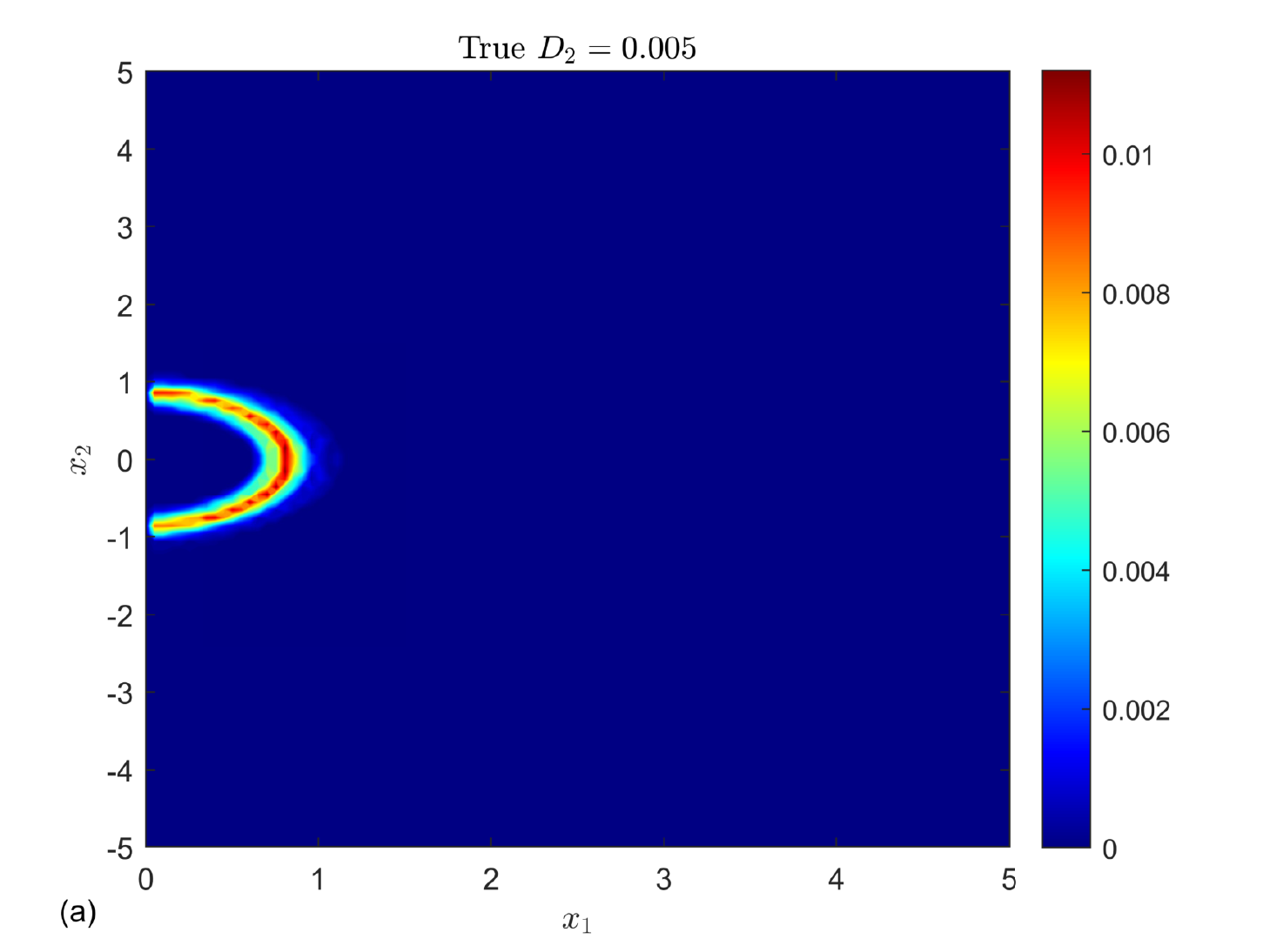}
    \includegraphics[width=5.3cm]{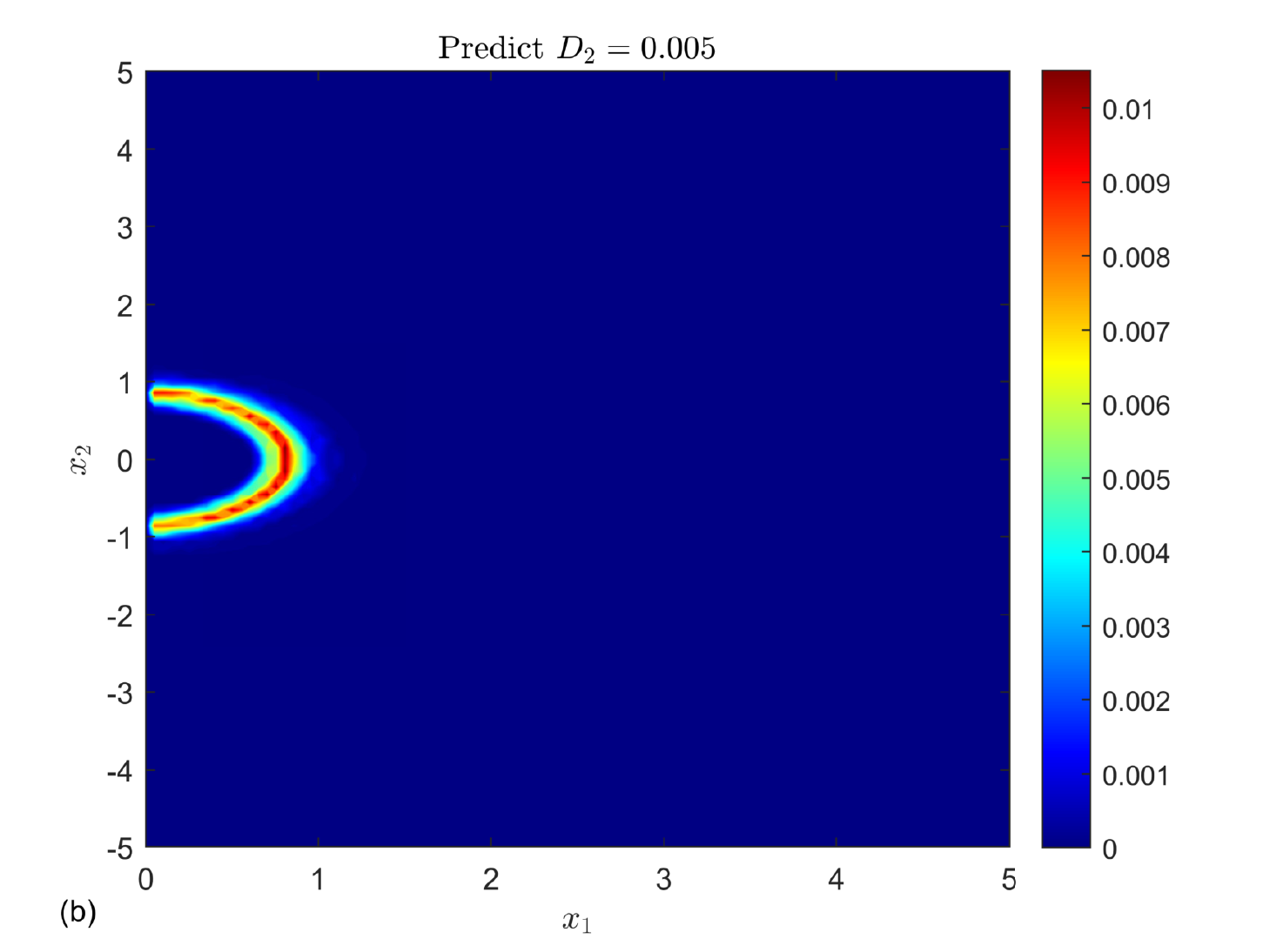}
    \includegraphics[width=5.3cm]{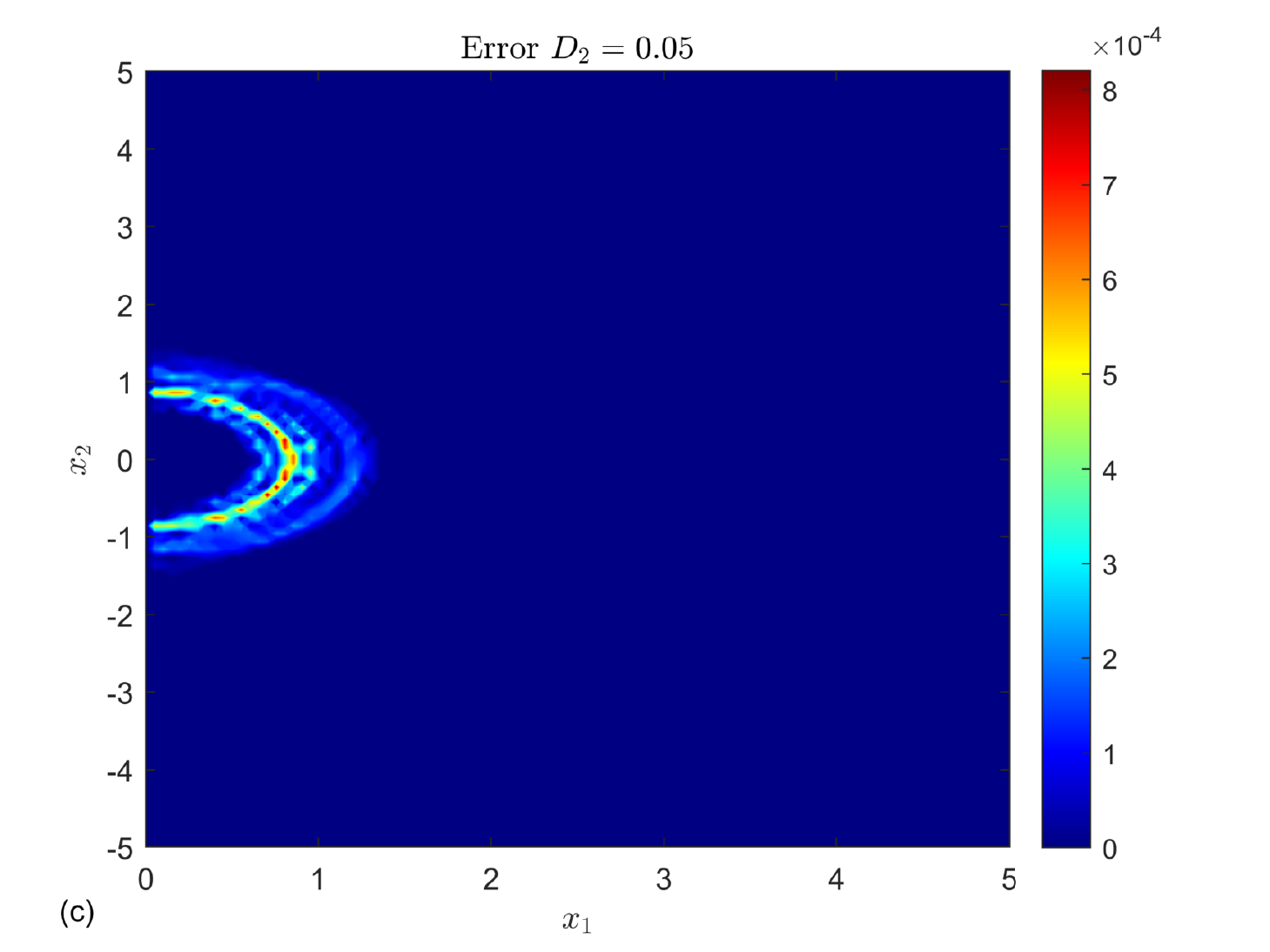} \\
    \includegraphics[width=5.3cm]{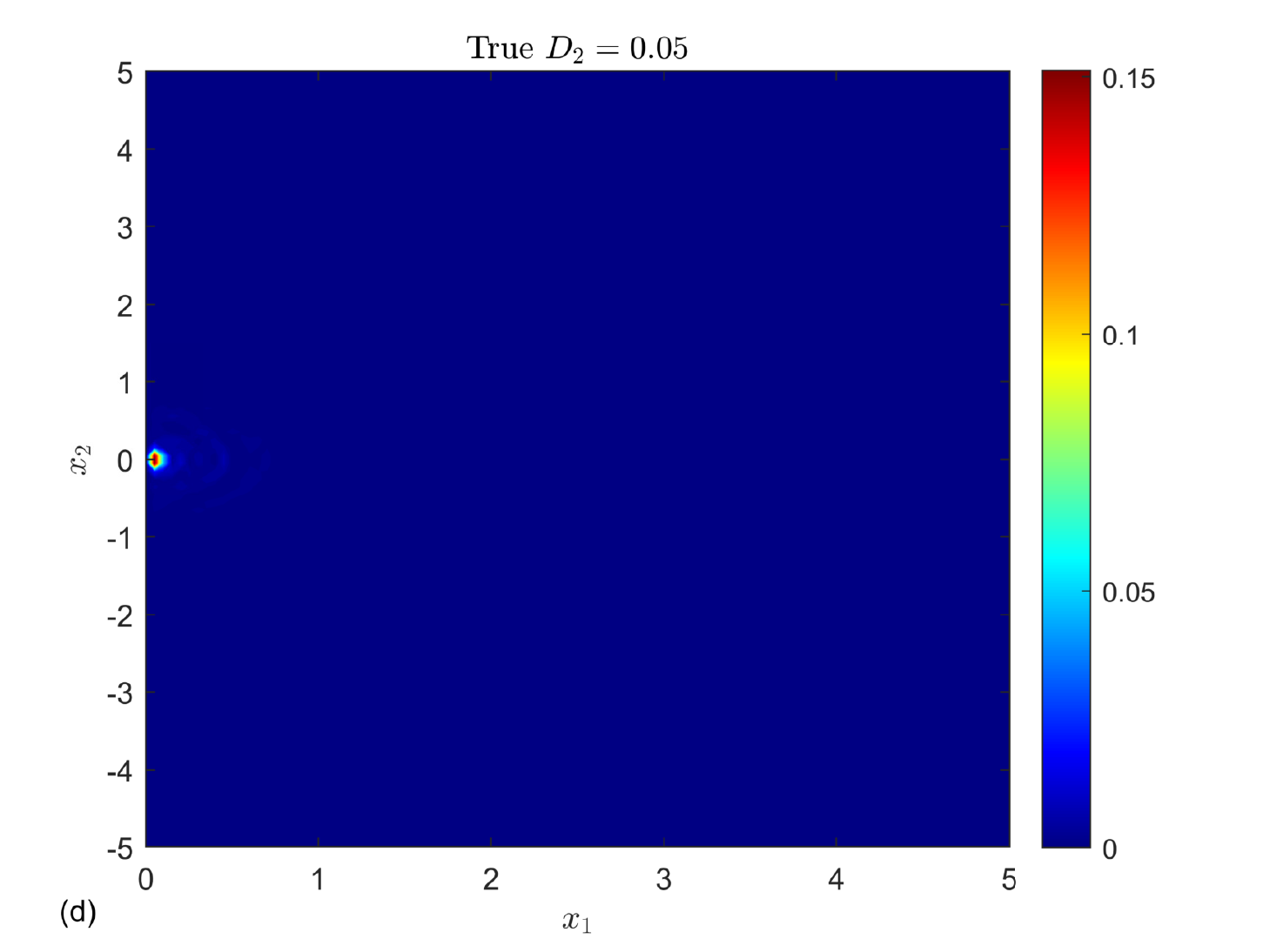}
    \includegraphics[width=5.3cm]{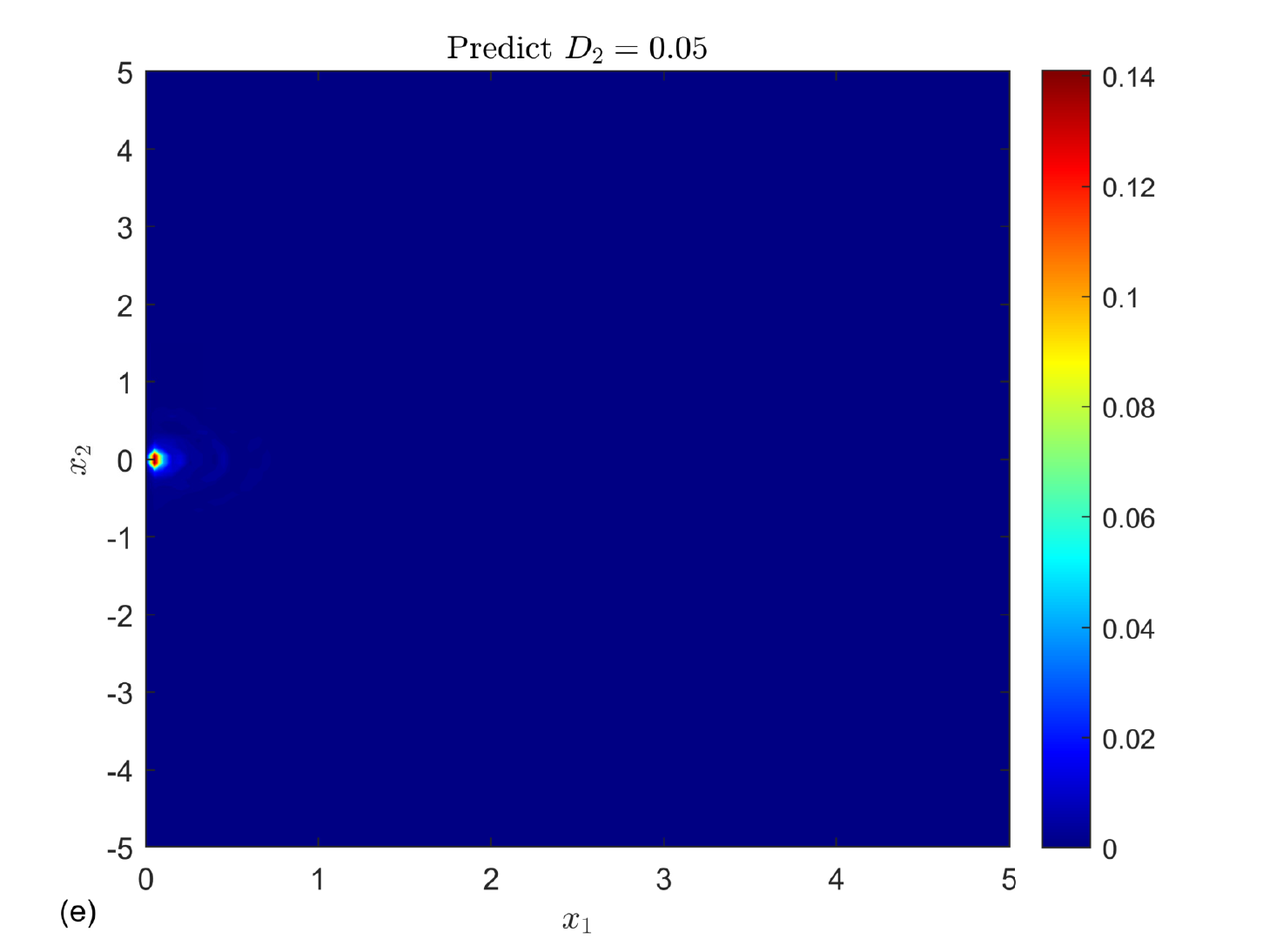}
    \includegraphics[width=5.3cm]{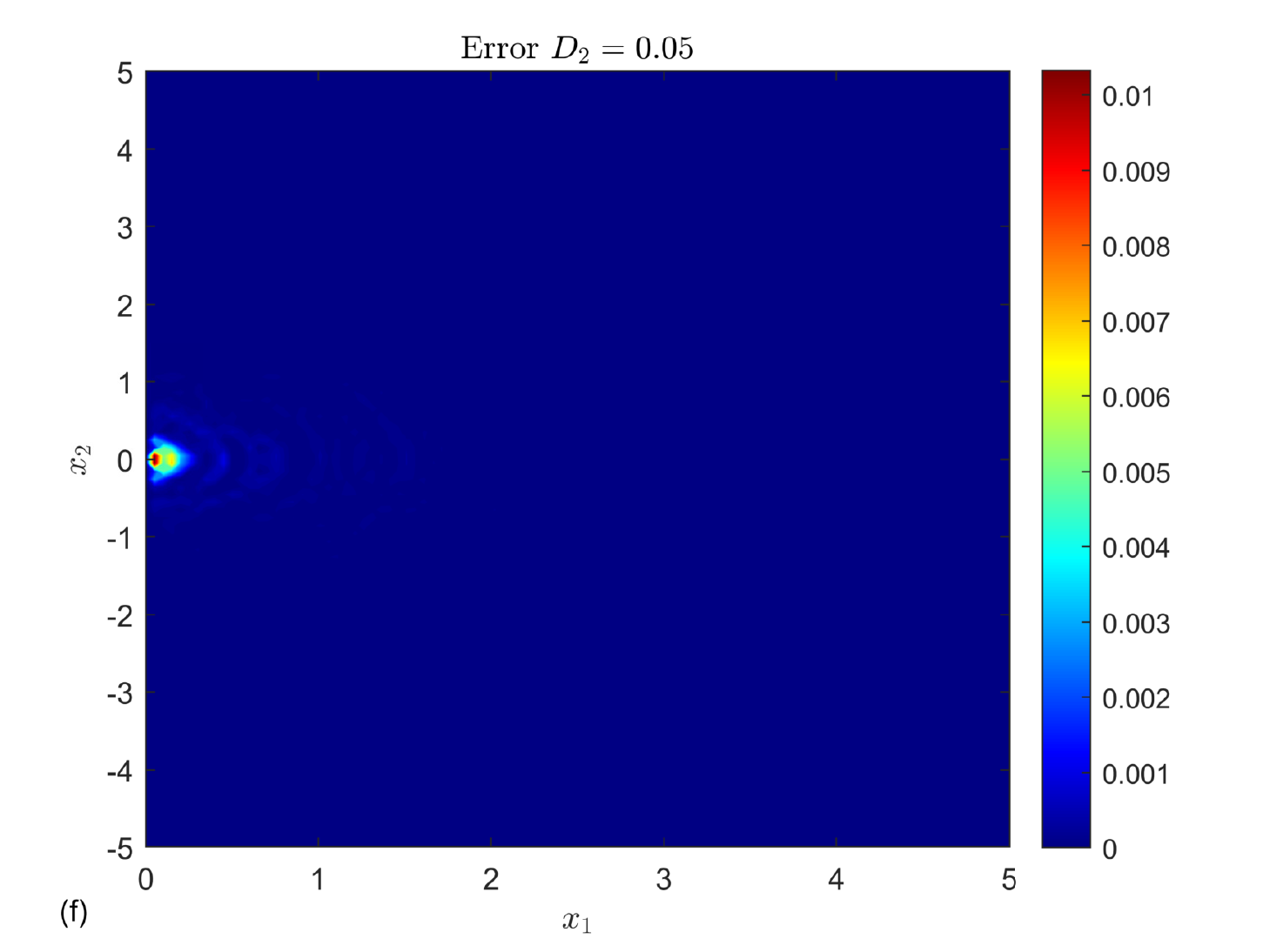}
    \caption{Predicting P-bifurcation Behavior of Rayleigh-Van der pol system Using Neural Network Functions.}
    \label{fig:6}
\end{figure}

In Fig.\ref{fig:6} we plot the top view of the predicted and true steady-state probability density functions. We employ Monte Carlo simulation to calculate the steady-state probability density function, which is solved iteratively using Eq.\eqref{28} and the neural network function, respectively. It can be seen that the noise intensity $D$ can significantly affect the shape of the stable most likely steady-state trajectory. As the noise intensity increases, the stable most probable steady-state trajectory collapses from a semicircle to a point, leading to the emergence of stochastic P bifurcations. In the prediction part, we use the neural network function to predict the stable most likely steady-state trajectory, and we can see that it has the same trend as the real image, and the occurrence of random P bifurcation is also observed. The two results are very close as shown in Fig.\ref{fig:6}(c,f). This shows that the stochastic differential equations fitted by the RBMS algorithm are effective, and the steady-state probability density function can be accurately calculated to observe the bifurcation behavior.

\section{Discussion}

The RBMS-II algorithm used in this paper has almost no hyperparameters, so we do not do sensitivity analysis of hyperparameters here. In the inverse problem of dynamics, we usually need to obtain a large amount of snapshot data of an unknown system, and realize fitting or prediction on this basis. In the actual engineering field, the collected data usually contain noise and disturbance, as shown in Eq.\eqref{30}. In order to eliminate and measure the impact of disturbance, a large number of scholars have developed a variety of methods, including perturbation method, bayesian method and uncertainty quantification method\cite{salahshour2021uncertain,he2003homotopy,beck1998updating, katafygiotis1998updating}. We try to mitigate the effects of noise and disturbances on the dynamical system inverse problem from the perspective of sampling methods.With the development of embedding physical information into neural networks and the popularization of big data technology, we believe this is a novel and important perspective.
\begin{equation}
    \dot{x}=f(t,x,\varepsilon ) ,x(t_{0} )=\eta (\varepsilon) \label{30}
\end{equation}

The method we use is highly dependent on data from 10,000 simulations starting from the same initial value. Just as the ancient Greek philosopher Heraclitus said that "one cannot step into the same river twice", in the field of practical engineering, the real initial values of the two simulations are at least slightly different. Therefore, we consider an initial value perturbation problem of the form Eq.\eqref{31}, combined with the properties of our method.

\begin{equation}
    \dot{x}=f(t,x,0) ,x(t_{0} )=\eta (\varepsilon) \label{31}
\end{equation}

\begin{figure}[htb]
    \centering
    \includegraphics[width=8.0cm]{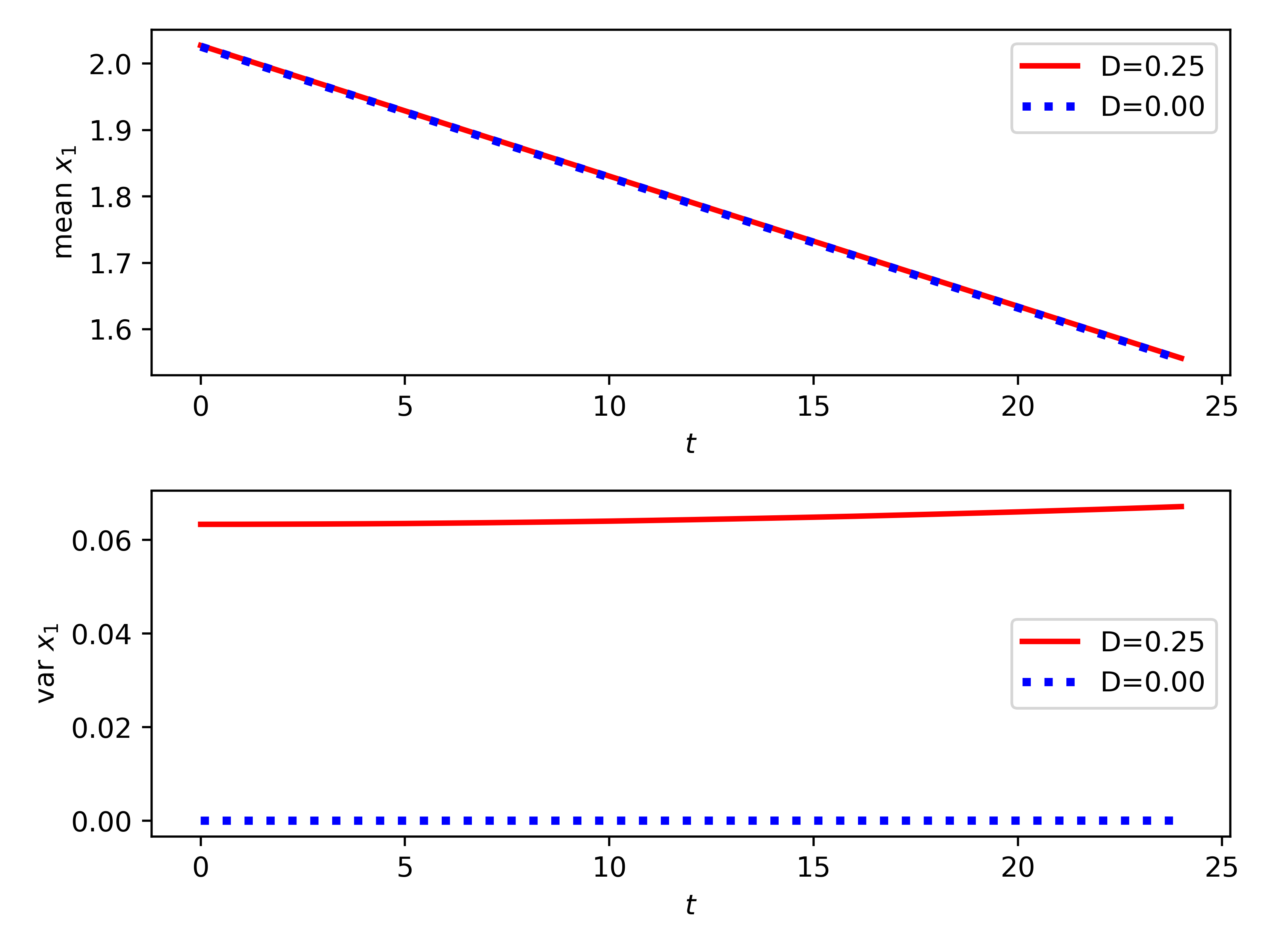} 
    \includegraphics[width=8.0cm]{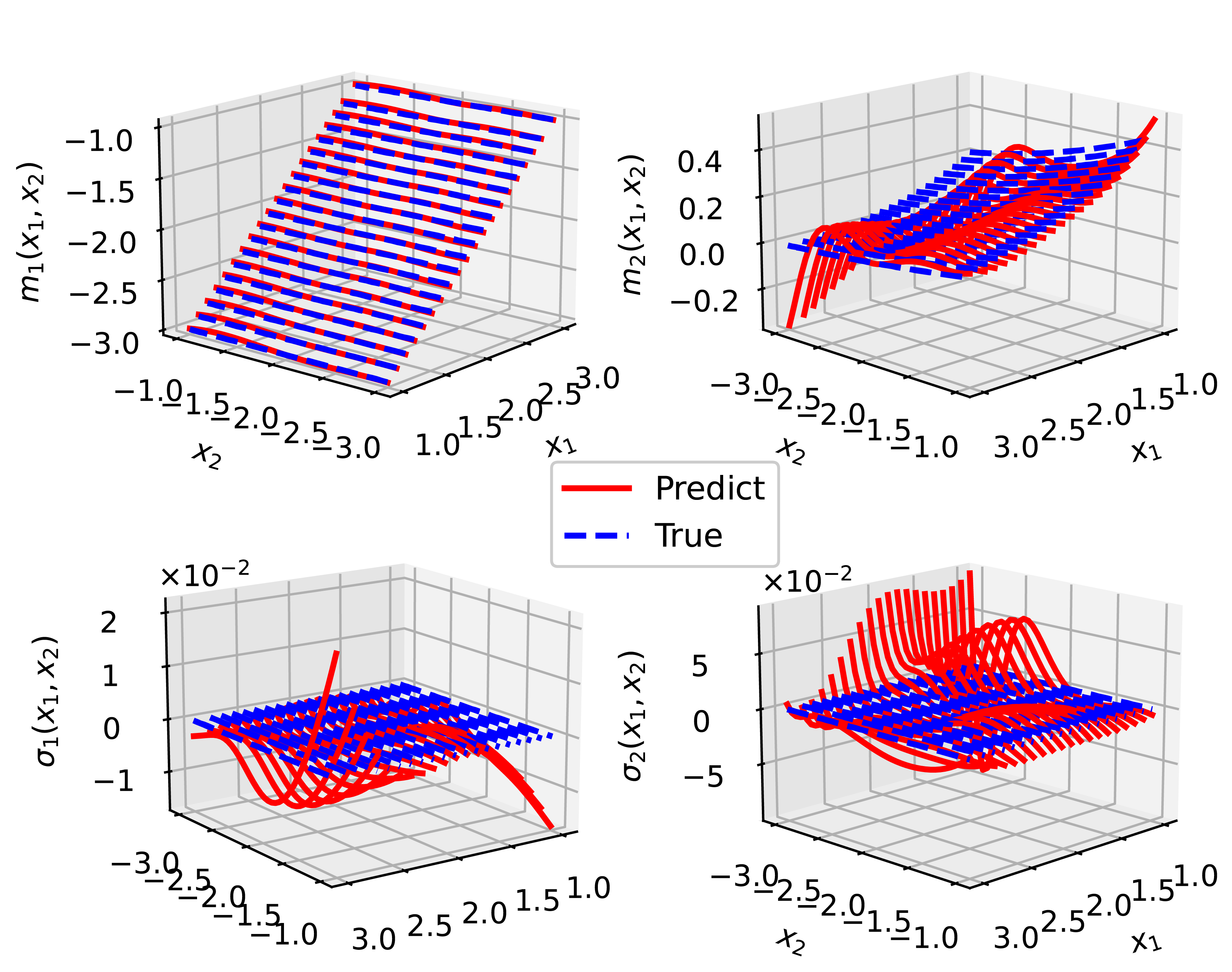} 
    
    \caption{Effect of initial value disturbance on Rayleigh-Van der pol system fitting task.}
    \label{fig:7}
\end{figure}

Let's take the Rayleigh-Van der pol system as an example, and set the noise at the initial value to obey the normal distribution $N(0,0.25^2)$. To evaluate the effect of perturbation at the initial value on the fitting task, we repeated the solution 10,000 times for an initial value with a time step of 25, and then observed the first and second moments of the solved data. In the subgraph on the left of Figure \ref{fig:7}, we compare the changes of the first-order moment and second-order moment with and without initial value disturbance. The noise we add follows a normal distribution with a mean of zero, so it has little effect on the first-order moments, but has a very large effect on the second-order moments. In order to further analyze the impact of the initial value disturbance on the fitting task, we use the noisy data to perform the fitting task of the Rayleigh-Van der pol system. In the right subplot of Figure \ref{fig:7}, we observe that both the drift term and the diffusion term are poorly fitted, especially the diffusion term. Therefore, we need to try to mitigate the bad effects of initial value perturbation.

\begin{figure}[htb]
    \centering
    \includegraphics[width=8.0cm]{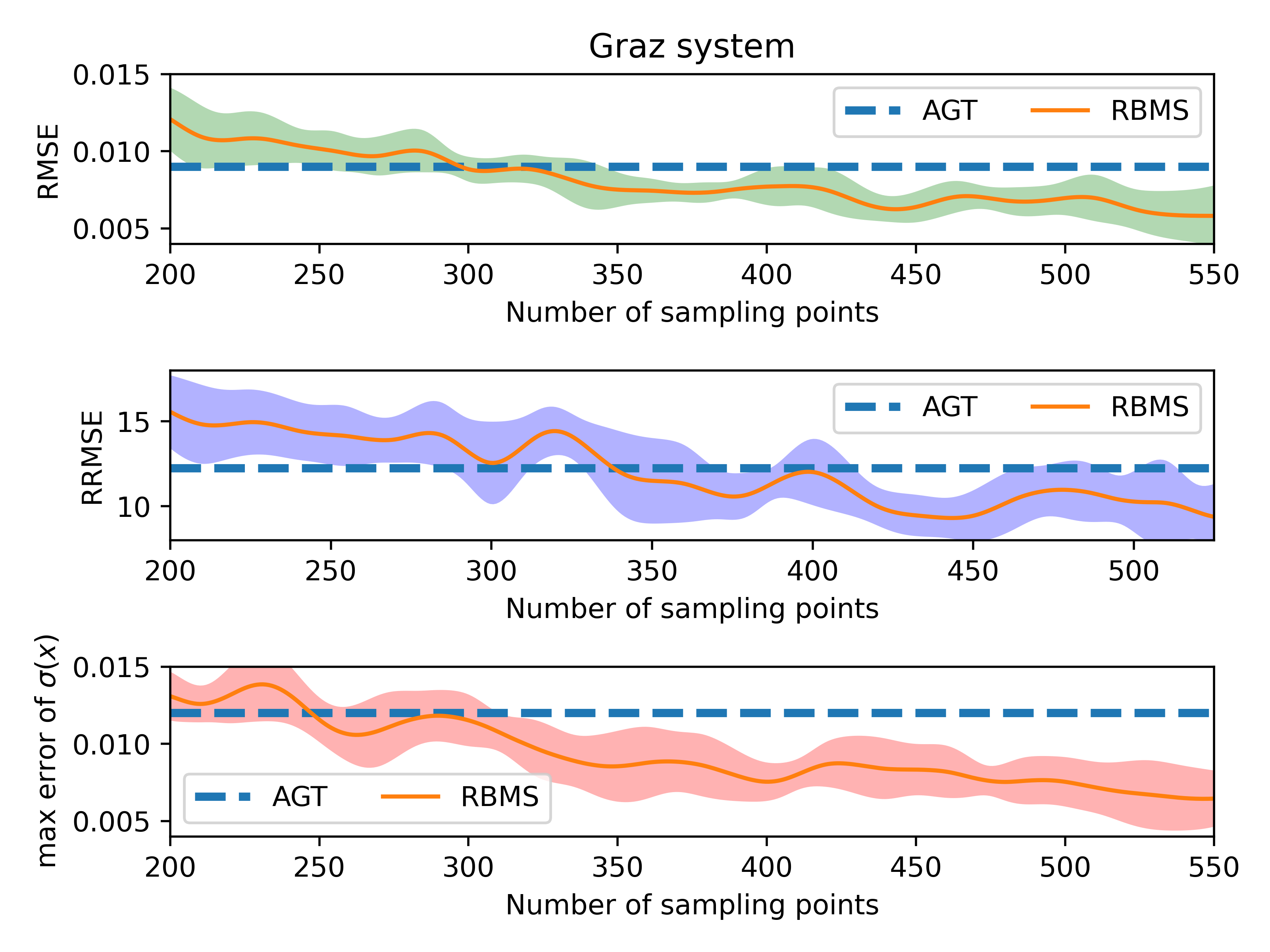}
    \includegraphics[width=8.0cm]{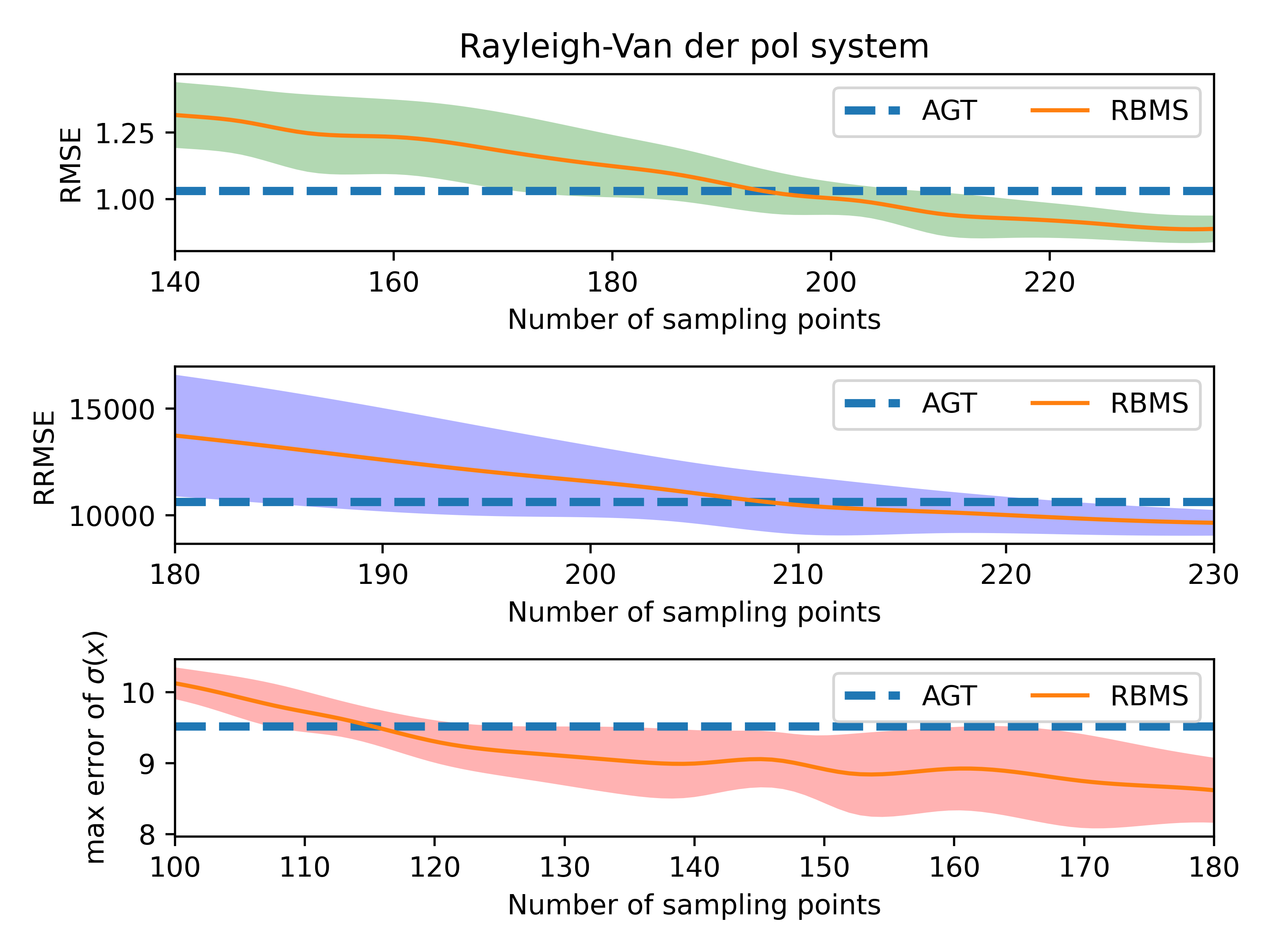} 
    \caption{Three kinds of error indicators for fitting Graz system and Rayleigh-Van der pol system using RBMS algorithm.}
    \label{fig:8}
\end{figure}

Combined with the above analysis of the impact of initial value disturbance, we establish three indicators to measure the fitting effect. In addition to using the RMSE and RRMSE to measure the overall error, we also added the maximum error of fitting the diffusion term as shown in Eq.\eqref{32}.

\begin{equation}
    max\ error\ of \ \sigma(X)=max\left \{ \left | \sigma ^{T}(X_{ij})- \theta  ^{P}(X_{ij}) \right |\right \} ,i=1,2,...,m;j=1,2,...,n.  \label{32}
\end{equation}

In Figure \ref{fig:8}, we show the fitting effects of the two systems using the RBMS algorithm under the monitoring of the three indicators.  AGT in the legend indicates that all points in the sampling grid (1600) are used to train the neural network, and  RBMS indicates that the neural network is trained using the sample points obtained by the RBMS algorithm. In order to reflect the stability of the method, for each system, we repeated the same experiment 10 times, and then drew the area graph with one standard deviation as the boundary.

It can be seen that in the fitting tasks of the grazing system and the Rayleigh-Van der pol system, using the RBMS algorithm can achieve better results than all point training when the sample size reaches about 400 and 220. Three indicators can reflect the advantages.

When the data contains noise, relying on too much data will lead to overfitting of the neural network relative to the noise, thereby significantly reducing the learning effect. Using a reasonable method to select a small number of samples for training can alleviate this overfitting, so the results produced by the RBMS algorithm are easy to understand and accept.

In the fitting task of the Rayleigh-Van der pol system, we observed very large fitting errors according to the RRMSE metric. This shows that compared with the grazing system, the second-order Rayleigh-Van der pol system system is more sensitive to initial value disturbance, and it is easier to fall into noise overfitting when using neural network fitting. Therefore, using only 220 sample points can achieve better results than all point training, and such results can be interpreted and accepted.

In summary, the RBMS algorithm we use has certain robustness in the task of fitting dynamical systems using deep learning methods. By using fewer samples, the impact of initial value disturbance on the fitting task can be alleviated to a certain extent.

\section{Conclusion}

More and more scholars try to apply the nonlinear system to the modeling of various problems because of its good characteristics. With the development of data science, data-driven modeling has become a powerful complementary method to first-principles modeling. The neural network embedded with physical information is the most popular data-driven modeling method, and there is still room for improvement in the part of the sampling algorithm. Based on the original RAR algorithm\cite{lu2021deepxde}, we propose a residual-based multimodal sampling algorithm (RBMS). After our proper tuning, the RBMS algorithm is a fully adaptive sampling algorithm with no hyperparameters.

We combined RBMS with SPINODE [23] to identify the physical information in two stochastic differential equations and then used a neural network as a ssurrogate function to analyse the stochastic dynamical behaviour of the system. The results show that our proposed RBMS algorithm requires only a small number of samples to achieve high accuracy, which means that the sampling process is efficient. In the discussion section, we test the robustness of the RBMS algorithm by adding initial value perturbations when solving for the simulated data. The results show that the RBMS algorithm can effectively mitigate the overfitting of the neural network to noise.

\section{Acknowledgements}
This work was supported by the National Natural Science
Foundation of China (NNSFC, Grant Nos. 12172291).

% 使用数据库文件 bib-example.bib   
% \bibliography{bib-example}

\begin{thebibliography}{10}

\bibitem{boechler2011bifurcation}
Neil Boechler, Georgios Theocharis, and C~Daraio.
\newblock Bifurcation-based acoustic switching and rectification.
\newblock {\em Nature materials}, 10(9):665--668, 2011.

\bibitem{antonio2012frequency}
Dario Antonio, Dami{\'a}n~H Zanette, and Daniel L{\'o}pez.
\newblock Frequency stabilization in nonlinear micromechanical oscillators.
\newblock {\em Nature communications}, 3(1):806, 2012.

\bibitem{strachan2013subharmonic}
B~Scott Strachan, Steven~W Shaw, and Oleg Kogan.
\newblock Subharmonic resonance cascades in a class of coupled resonators.
\newblock {\em Journal of Computational and Nonlinear Dynamics}, 8(4), 2013.

\bibitem{kecik2014parametric}
Krzysztof Kecik and Marcin Kapitaniak.
\newblock Parametric analysis of magnetorheologically damped pendulum vibration
  absorber.
\newblock {\em International Journal of Structural Stability and Dynamics},
  14(08):1440015, 2014.

\bibitem{vakakis2008nonlinear}
Alexander~F Vakakis, Oleg~V Gendelman, Lawrence~A Bergman, D~Michael McFarland,
  Ga{\"e}tan Kerschen, and Young~Sup Lee.
\newblock {\em Nonlinear targeted energy transfer in mechanical and structural
  systems}, volume 156.
\newblock Springer Science \& Business Media, 2008.

\bibitem{green2012benefits}
PL~Green, K~Worden, K~Atallah, and ND~Sims.
\newblock The benefits of duffing-type nonlinearities and electrical
  optimisation of a mono-stable energy harvester under white gaussian
  excitations.
\newblock {\em Journal of sound and vibration}, 331(20):4504--4517, 2012.

\bibitem{amin2012powering}
M~Amin~Karami and Daniel~J Inman.
\newblock Powering pacemakers from heartbeat vibrations using linear and
  nonlinear energy harvesters.
\newblock {\em Applied Physics Letters}, 100(4):042901, 2012.

\bibitem{quinn2011energy}
D~Dane Quinn, Angela~L Triplett, Alexander~F Vakakis, and Lawrence~A Bergman.
\newblock Energy harvesting from impulsive loads using intentional essential
  nonlinearities.
\newblock {\em Journal of Vibration and Acoustics}, 133(1), 2011.

\bibitem{psichogios1992hybrid}
Dimitris~C Psichogios and Lyle~H Ungar.
\newblock A hybrid neural network-first principles approach to process
  modeling.
\newblock {\em AIChE Journal}, 38(10):1499--1511, 1992.

\bibitem{rico1992discrete}
Ramiro Rico-Martinez, K~Krischer, IG~Kevrekidis, MC~Kube, and JL~Hudson.
\newblock Discrete-vs. continuous-time nonlinear signal processing of cu
  electrodissolution data.
\newblock {\em Chemical Engineering Communications}, 118(1):25--48, 1992.

\bibitem{thompson1994modeling}
Michael~L Thompson and Mark~A Kramer.
\newblock Modeling chemical processes using prior knowledge and neural
  networks.
\newblock {\em AIChE Journal}, 40(8):1328--1340, 1994.

\bibitem{liu2022complex}
Qi~Liu, Yong Xu, J{\"u}rgen Kurths, and Xiaochuan Liu.
\newblock Complex nonlinear dynamics and vibration suppression of conceptual
  airfoil models: A state-of-the-art overview.
\newblock {\em Chaos: An Interdisciplinary Journal of Nonlinear Science},
  32(6):062101, 2022.

\bibitem{noel2017nonlinear}
Jean-Philippe No{\"e}l and Ga{\"e}tan Kerschen.
\newblock Nonlinear system identification in structural dynamics: 10 more years
  of progress.
\newblock {\em Mechanical Systems and Signal Processing}, 83:2--35, 2017.

\bibitem{pei2015demonstration}
Jin-Song Pei and Sami~F Masri.
\newblock Demonstration and validation of constructive initialization method
  for neural networks to approximate nonlinear functions in engineering
  mechanics applications.
\newblock {\em Nonlinear Dynamics}, 79:2099--2119, 2015.

\bibitem{ayala2016cascaded}
Helon Vicente~Hultmann Ayala and Leandro dos Santos~Coelho.
\newblock Cascaded evolutionary algorithm for nonlinear system identification
  based on correlation functions and radial basis functions neural networks.
\newblock {\em Mechanical Systems and Signal Processing}, 68:378--393, 2016.

\bibitem{tavakolpour2015parametric}
AR~Tavakolpour-Saleh, SAR Nasib, A~Sepasyan, and SM~Hashemi.
\newblock Parametric and nonparametric system identification of an experimental
  turbojet engine.
\newblock {\em Aerospace Science and Technology}, 43:21--29, 2015.

\bibitem{paduart2010identification}
Johan Paduart, Lieve Lauwers, Jan Swevers, Kris Smolders, Johan Schoukens, and
  Rik Pintelon.
\newblock Identification of nonlinear systems using polynomial nonlinear state
  space models.
\newblock {\em Automatica}, 46(4):647--656, 2010.

\bibitem{dreesen2015decoupling}
Philippe Dreesen, Mariya Ishteva, and Johan Schoukens.
\newblock Decoupling multivariate polynomials using first-order information and
  tensor decompositions.
\newblock {\em SIAM Journal on Matrix Analysis and Applications},
  36(2):864--879, 2015.

\bibitem{greydanus2019hamiltonian}
Samuel Greydanus, Misko Dzamba, and Jason Yosinski.
\newblock Hamiltonian neural networks.
\newblock {\em Advances in neural information processing systems}, 32, 2019.

\bibitem{cranmer2020lagrangian}
Miles Cranmer, Sam Greydanus, Stephan Hoyer, Peter Battaglia, David Spergel,
  and Shirley Ho.
\newblock Lagrangian neural networks.
\newblock {\em arXiv preprint arXiv:2003.04630}, 2020.

\bibitem{raissi2019physics}
Maziar Raissi, Paris Perdikaris, and George~E Karniadakis.
\newblock Physics-informed neural networks: A deep learning framework for
  solving forward and inverse problems involving nonlinear partial differential
  equations.
\newblock {\em Journal of Computational physics}, 378:686--707, 2019.

\bibitem{chen2018neural}
Ricky~TQ Chen, Yulia Rubanova, Jesse Bettencourt, and David~K Duvenaud.
\newblock Neural ordinary differential equations.
\newblock {\em Advances in neural information processing systems}, 31, 2018.

\bibitem{o2022stochastic}
Jared O'Leary, Joel~A Paulson, and Ali Mesbah.
\newblock Stochastic physics-informed neural ordinary differential equations.
\newblock {\em Journal of Computational Physics}, 468:111466, 2022.

\bibitem{lu2021deepxde}
Lu~Lu, Xuhui Meng, Zhiping Mao, and George~Em Karniadakis.
\newblock Deepxde: A deep learning library for solving differential equations.
\newblock {\em SIAM review}, 63(1):208--228, 2021.

\bibitem{nabian2021efficient}
Mohammad~Amin Nabian, Rini~Jasmine Gladstone, and Hadi Meidani.
\newblock Efficient training of physics-informed neural networks via importance
  sampling.
\newblock {\em Computer-Aided Civil and Infrastructure Engineering},
  36(8):962--977, 2021.

\bibitem{wu2023comprehensive}
Chenxi Wu, Min Zhu, Qinyang Tan, Yadhu Kartha, and Lu~Lu.
\newblock A comprehensive study of non-adaptive and residual-based adaptive
  sampling for physics-informed neural networks.
\newblock {\em Computer Methods in Applied Mechanics and Engineering},
  403:115671, 2023.

\bibitem{gao2023active}
Wenhan Gao and Chunmei Wang.
\newblock Active learning based sampling for high-dimensional nonlinear partial
  differential equations.
\newblock {\em Journal of Computational Physics}, 475:111848, 2023.

\bibitem{tang2023pinns}
Kejun Tang, Xiaoliang Wan, and Chao Yang.
\newblock Das-pinns: A deep adaptive sampling method for solving
  high-dimensional partial differential equations.
\newblock {\em Journal of Computational Physics}, 476:111868, 2023.

\bibitem{hanna2022residual}
John~M Hanna, Jose~V Aguado, Sebastien Comas-Cardona, Ramzi Askri, and Domenico
  Borzacchiello.
\newblock Residual-based adaptivity for two-phase flow simulation in porous
  media using physics-informed neural networks.
\newblock {\em Computer Methods in Applied Mechanics and Engineering},
  396:115100, 2022.

\bibitem{zeng2022adaptive}
Shaojie Zeng, Zong Zhang, and Qingsong Zou.
\newblock Adaptive deep neural networks methods for high-dimensional partial
  differential equations.
\newblock {\em Journal of Computational Physics}, 463:111232, 2022.

\bibitem{peng2022rang}
Wei Peng, Weien Zhou, Xiaoya Zhang, Wen Yao, and Zheliang Liu.
\newblock Rang: A residual-based adaptive node generation method for
  physics-informed neural networks.
\newblock {\em arXiv preprint arXiv:2205.01051}, 2022.

\bibitem{zapf2022investigating}
Bastian Zapf, Johannes Haubner, Miroslav Kuchta, Geir Ringstad, Per~Kristian
  Eide, and Kent-Andre Mardal.
\newblock Investigating molecular transport in the human brain from mri with
  physics-informed neural networks.
\newblock {\em Scientific Reports}, 12(1):15475, 2022.

\bibitem{julier1997new}
Simon~J Julier and Jeffrey~K Uhlmann.
\newblock New extension of the kalman filter to nonlinear systems.
\newblock In {\em Signal processing, sensor fusion, and target recognition VI},
  volume 3068, pages 182--193. Spie, 1997.

\bibitem{noy1975stability}
Imanuel Noy-Meir.
\newblock Stability of grazing systems: an application of predator-prey graphs.
\newblock {\em The Journal of Ecology}, pages 459--481, 1975.

\bibitem{zhang2019noise}
Hongxia Zhang, Wei Xu, Youming Lei, and Yan Qiao.
\newblock Noise-induced vegetation transitions in the grazing ecosystem.
\newblock {\em Applied Mathematical Modelling}, 76:225--237, 2019.

\bibitem{eugene1965relation}
Wong Eugene and Zakai Moshe.
\newblock On the relation between ordinary and stochastic differential
  equations.
\newblock {\em International Journal of Engineering Science}, 3(2):213--229,
  1965.

\bibitem{duan2015introduction}
Jinqiao Duan.
\newblock {\em An introduction to stochastic dynamics}, volume~51.
\newblock Cambridge University Press, 2015.

\bibitem{paget1937vibration}
AL~Paget.
\newblock Vibration in steam turbine buckets and damping by impacts.
\newblock {\em Engineering}, 143:305--307, 1937.

\bibitem{foale1994bifurcations}
S~Foale and SR~Bishop.
\newblock Bifurcations in impact oscillations.
\newblock {\em Nonlinear dynamics}, 6:285--299, 1994.

\bibitem{ibrahim2009vibro}
Raouf~A Ibrahim.
\newblock {\em Vibro-impact dynamics: modeling, mapping and applications},
  volume~43.
\newblock Springer Science \& Business Media, 2009.

\bibitem{salahshour2021uncertain}
Soheil Salahshour, Ali Ahmadian, Bruno~A Pansera, and Massimiliano Ferrara.
\newblock Uncertain inverse problem for fractional dynamical systems using
  perturbed collage theorem.
\newblock {\em Communications in Nonlinear Science and Numerical Simulation},
  94:105553, 2021.

\bibitem{he2003homotopy}
Ji-Huan He.
\newblock Homotopy perturbation method: a new nonlinear analytical technique.
\newblock {\em Applied Mathematics and computation}, 135(1):73--79, 2003.

\bibitem{beck1998updating}
James~L Beck and Lambros~S Katafygiotis.
\newblock Updating models and their uncertainties. i: Bayesian statistical
  framework.
\newblock {\em Journal of Engineering Mechanics}, 124(4):455--461, 1998.

\bibitem{katafygiotis1998updating}
Lambros~S Katafygiotis and Jim~L Beck.
\newblock Updating models and their uncertainties. ii: Model identifiability.
\newblock {\em Journal of Engineering Mechanics}, 124(4):463--467, 1998.

\end{thebibliography}

\end{document}